\documentclass[a4paper, reqno]{amsart}
\usepackage{amssymb, amsfonts, amsthm, amsmath, mathtext, latexsym,graphicx}
\usepackage{color}
\newcommand{\dsize}{\textstyle}
\newcommand{\U}{\overset{\circ}U}
\newcommand{\R}{{\mathbb R}}

\newcommand{\Z}{{\mathbb Z}}
\newcommand{\N}{{\mathbb N}}

\newcommand{\C}{{\mathbb C}}
\newcommand{\D}{\displaystyle}  

\renewcommand{\mod}{{\rm mod}\,}
\newcommand{\re}{{\rm Re}\,}
\newcommand{\im}{{\rm Im}\,}

\newcommand{\ch}{{\mathrm cosh}\,}
\newcommand{\sh}{{\mathrm sinh}\,}

\theoremstyle{plain}
\newtheorem{Th}{Theorem}[section]
\newtheorem{Le}{Lemma}[section]
\newtheorem{Pro}{Proposition}[section]
\newtheorem{Cor}{Corollary}[section]

\theoremstyle{definition}
\newtheorem{Rem}{Remark}[section]

\numberwithin{equation}{section}
\title{The complex WKB method for difference equations and Airy functions}
\author{Alexander Fedotov and Fr{\'e}d{\'e}ric Klopp}
\address[Alexander Fedotov]{St. Petersburg State University, 
7/9 Universitetskaya nab., St.Petersburg, 199034, Russia}
\email{
{a.fedotov@spbu.ru}}
\address[Fr{\'e}d{\'e}ric Klopp]{Sorbonne Universit{\'e}, Universit{\'e} Paris Diderot,
  CNRS, Institut de Math{\'e}matiques de Jussieu - Paris Rive Gauche ,
  F-75005, Paris, France}
\email{
{frederic.klopp@imj-prg.fr}}
\begin{document}
\begin{abstract}
  We consider the difference Schr{\"o}dinger equation
  $\psi(z+h)+\psi(z-h)+ v(z)\psi(z)=0$ where $z$ is a complex
  variable, $h>0$ is a parameter, and $v$ is an analytic function. As
  $h\to0$ analytic solutions to this equation have a standard
  quasiclassical behavior near the points where $v(z)\ne\pm 2$. We
  study analytic solutions near the points $z_0$ satisfying
  $v(z_0)=\pm2$ and $v'(z_0)\ne 0$. For the finite difference
  equation, these points are the natural analogues of the simple
  turning points defined for the differential equation
  $-\psi''(z)+v(z)\psi(z)=0$. In an $h$-independent neighborhood of
  such a point, we derive uniform asymptotic expansions for analytic
  solutions to the difference equation.
\end{abstract}
\thanks{The work was supported by CNRS, France, and the Russian
  Foundation for Basic Research under the grant No 17-51-150008.}
\maketitle
\section{Introduction, Preliminaries, and Main Results}
\label{sec:introduction}
\subsection{The problem}
We study analytic solutions to the difference Schr{\"o}dinger
equation
\begin{equation}
  \label{eq:main}
  \psi(z+h)+\psi(z-h)+ v(z)\psi(z)=0
\end{equation}
where $z$ is a complex variable, $h$ is a positive parameter and $v$
is an analytic function. We describe their
asymptotics as $h\to 0$.\\
Note that the parameter $h$ is a standard quasiclassical
parameter. Indeed, formally,
$\psi(z+h)=\sum_{l=0}^\infty \frac{h^l}{l!}
\frac{d^l\psi}{dz^l}(z)=e^{h\frac{d}{dz}}\psi(z)$, and $h$ can be
regarded as a small parameter in front of the
derivative.\\
One encounters difference equations in the complex plane in many
fields of mathematics and physics.  For example, they arise when
studying an electron in a crystal submitted to a constant magnetic
field (e.g.,~\cite{GHT:89}), wave scattering by wedges
(e.g.,~\cite{B-L-G:2008}) and one-dimensional quasi-periodic
differential Schr{\"o}dinger equations with two frequencies
(e.g.,~\cite{F-K:02}). The quasiclassical case corresponds
respectively to the cases of a small magnetic field, of a thin wedge
and the case where one frequency is small with respect to another.\\
The quasiclassical theory of difference equations in the complex plane
can also be useful to study orthogonal polynomials,
see section~\ref{rel-res}.\\
The quasiclassical asymptotics of analytic solutions to ordinary
differential equations in the complex plane are described by the
well-known complex WKB method (see, e.g.,~\cite{W:87, Fe:93}). The
complex WKB method for difference equations was developed
in~\cite{B-F:94, F-Shch:15, F-Shch:18}.\\
The present paper is devoted to uniform asymptotic formulas describing
analytic solutions to~\eqref{eq:main} in $h$-independent complex
neighborhoods of simple turning points (see
sections~\ref{sec:complex-momentum} and~\ref{sec:compl-moment-near}).
The results of this paper were partially announced in~\cite{F-K:18}.

The idea to study the asymptotics of solutions to a difference 
equation in a complex neighborhood of a turning point appears to be very natural.
One can say that this idea and the techniques developed to get the asymptotics
are the main analytic innovations of the paper.

In the next sections, we first recall some basic definitions and
statements of the complex WKB method for difference equations. Next,
we introduce a few objects needed to formulate our results that we
then state.\\
We  assume that $v$ is analytic on a disk $U\subset\C$.\\
Below, a neighborhood is a $\delta$-neighborhood, in particular, a
neighborhood of a point is an open disk with the center at this point.
\subsection{A very short introduction to the complex WKB method}
\label{sec:very-short-intr}
Here, following~\cite{F-Shch:15, F-Shch:18}, we briefly describe basic 
definitions and results of the complex WKB method for difference equations.
\subsubsection{The complex momentum}
\label{sec:complex-momentum}
The main analytic object of the complex WKB method is the {\it complex
  momentum} $p$. For~\eqref{eq:main} it is defined by the formula
\begin{equation}
  \label{eq:momentum}
  2\cos p+v(z)=0.
\end{equation}
It is a multivalued analytic function on $U$. At its branching points
$\cos p(z)\in\{\pm1\}$, thus, $v(z)\in\{\pm 2\}$.\\
In analogy with the glossary of the complex WKB method for
differential equations, the points where $v(z)\in\{\pm2\} $
are called {\it turning points}.\\
A set $D\subset U$ is {\it regular} if $v(z)\ne \pm 2$ in $D$.
\subsubsection{The main theorem of the complex WKB method}
\label{sec:main-theorem-complex}
As in the case of differential equations, one of the main geometric notions of
the complex WKB method is the {\it canonical domain}. In this paper 
we do not use it directly, and the reader needs to keep in mind only that the 
canonical domains are regular, simply connected domains independent of $h$, 
and that one has the following two theorems.
\begin{Th}
  \label{th:existence-can-dom}
  Any regular point belongs to a canonical domain.
\end{Th}
\noindent The proof of this statement repeats the proof of Lemma 5.2
from~\cite{F-K:05}.
\begin{Th}
  Let $K\subset U$ be a canonical domain, let $z_0\in K$, and let $p$ be a
  branch of the complex momentum analytic in $K$. Then there exist
  solutions $\psi_\pm$ to~\eqref{eq:main} analytic in $K$ and such
  that as $h\to 0$
  \begin{equation}
    \label{eq:main:wkb}
    \psi_\pm(z)=\frac1{\sqrt{\sin(p(z))}}\, e^{\pm\frac{i}{h}
      \int_{z_0}^zp(z)\,dz+o(1)}.
    \quad z\in K.
  \end{equation}
  This asymptotic representation is locally uniform in $K$.
\end{Th}
\noindent In \cite{F-Shch:15} this theorem was proved for $v$ analytic in 
bounded domains, and in \cite{F-Shch:18} it was proved in the case
where $v$ is a trigonometric polynomial.\\
Note that, by definition, at a turning point of $p$, one has
$\sin p(z)=0$. Thus, representation~\eqref{eq:main:wkb} cannot be
valid in a neighborhood of a turning point.
\begin{Rem}
  \label{rem:1}
  For the differential equation $-h^2\psi''(z)+v(z)\psi(z)=0$,
  formula~\eqref{eq:main:wkb} has to be replaced with (see, e.g.~\cite{W:87})
  $ \psi_\pm(z)=\frac1{\sqrt{p(z)}}\,
  e^{\pm\frac{i}{h}\int_{z_0}^zp(z)\,dz+o(1)}$,
  where the complex momentum is defined by the relation $p^2+v(z)=0$,
  i.e., as for~\eqref{eq:main}, by the symbol of the equation.
\end{Rem}
\subsection{The complex momentum and the conformal mapping $\zeta$}
\label{ss:zeta}
Here, we discuss properties of the complex momentum that we use throughout this paper.
These properties easily follow from the definition of $p$.
\subsubsection{Analytic branches of the complex momentum}
\label{sec:analyt-branch-compl}
Let $p_0$ be a branch of the complex momentum analytic in a regular
simply connected domain $D$. Then, an analytic function
$\tilde p:\ D\to\C$ is a branch of the complex momentum if and only if
there exists $s\in\{\pm1\}$ and $n\in \Z$ such that
\begin{equation}
  \label{eq:p-branches}
  \tilde p(z)=s p_0(z)+2\pi n,\quad \forall z\in D. 
\end{equation}
\subsubsection{The values of $p$ at turning points}
\label{sec:values-at-branch}
We note that $z_0\in U$ is a turning point if and only if
$p(z_0)=0\,\mod \pi$.  A simple transformation of the equation shows
that it suffices to consider the case where $p(z_0)=0\,\mod
2\pi$. Indeed, for $\psi$, a solution to~\eqref{eq:main}, we set
$\phi(z)=e^{i\pi z/h}\psi(z)$. Then, $\phi$ satisfies equation
\begin{equation}
  \label{eq:aux}
  \phi(z+h)+\phi(z-h)- v(z)\phi(z)=0.
\end{equation}
The complex momenta for equations~\eqref{eq:main} and~\eqref{eq:aux}
differ by $\pi\, \mod 2\pi$, and $z_0$ is a turning point
for~\eqref{eq:main} if and only if it is a turning point
for~\eqref{eq:aux}.
\subsubsection{The complex momentum near a turning point}
\label{sec:compl-moment-near}
Let $z_0$ be a turning point. If $v'(z_0)\ne0$, we call the turning
point $z_0$ {\it simple}. In this case, the complex momentum is
analytic in $\tau=\sqrt{z-z_0}$ in a neighborhood of $0$, and as
$\tau\to0$ any of its analytic branches admits a representation of the
form
\begin{equation}
  \label{eq:p-z}
  p(z)=p(z_0)+k_1\tau+O(\tau^2),\quad \tau=\sqrt{z-z_0}, \qquad k_1\ne0.
\end{equation}
\subsubsection{Our assumptions}
\label{assumpt:p}
From now on, we assume that
\begin{itemize}
\item in the disk $U$, there exists a single turning point, namely its
  center $z_0$, and it is simple;
\item $p(z_0)=0\,\mod 2\pi$.
\end{itemize}
\subsubsection{The function $\zeta$}
\label{sec:function-zeta}
The function $\zeta$ we describe here plays an important role in the asymptotic analysis
of~\eqref{eq:main} near turning points.\\
We cut $U$ from $z_0$  to a point of its boundary along a simple
curve  and denote the thus obtained domain by $U'$. In $U'$, we fix an analytic branch 
$p$ of the complex momentum.

We have $p(z_0)=2\pi n$, \ $n\in\Z$. Clearly, $p(z)-2\pi n$ also is a  
branch of the complex momentum analytic in $U'$.
So, we can and do assume that $p(z_0)=0$.

Let us fix in $U'$  an analytic branch $\zeta$ of the function
\begin{equation}
  \label{zeta}
  z\mapsto
  \left(\frac{3}{2i}\int_{z_0}^zp(z)\,dz\right)^{\frac23}.
\end{equation}
This branch  is actually analytic in $U$. One has
$\zeta(z_0)=0$, and $\zeta'(z_0)\ne 0$. \\ 
\begin{Rem} 
There are three different analytic branches of function~\eqref{zeta}: 
they equal $e^{4\pi i j /3}\zeta$, $j\in\Z_3=\Z/3\Z$. 
The set of these branches is independent of the curve along which we  
cut $U$ to get $U'$ and on the precise choice of the branch $p$.
\end{Rem}
We note that the definition of $\zeta$  implies that it satisfies one of the
two equations
\begin{equation}\label{dif-eq:zeta}
 \sqrt{\zeta(z)}\zeta'(z)=\pm ip(z), \qquad z\in U.
\end{equation}
Possibly reducing $U$ somewhat, we can and do assume that
\begin{itemize}
\item $\zeta$ is a bi-analytic bijection of $U$ onto its image.
\end{itemize}
\subsection{Basic facts on Airy functions}
\label{BFAF}
\begin{figure}
  \centering
  \includegraphics[height=4cm]{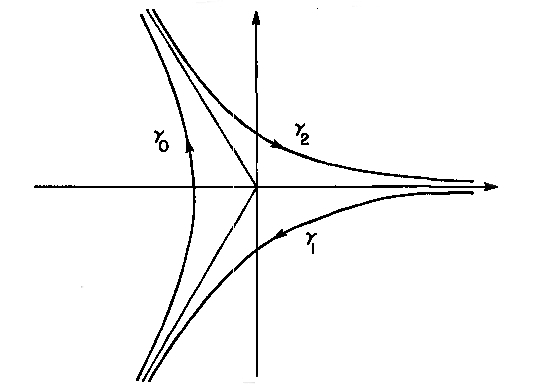}
  \caption{Integration paths}\label{fig:Airy}
\end{figure}

The equation
\begin{equation}\label{eq:Airy}
  w''(\zeta)=\zeta w(\zeta), \qquad \zeta\in\C,
\end{equation}
is the Airy equation. Its solutions are {\it Airy functions}.\\
Let $(\gamma_j)_{j\in\Z_3}$ be the curves shown in
Fig.~\ref{fig:Airy} borrowed from~\cite{W}; $\gamma_0$ is asymptotic to the half-lines
$e^{\pm 2i\pi/3}\R_+$, \ $\R_+=[0,+\infty)\subset\R$; for $j\in\Z_3$,
rotating $\gamma_0$ around $0$ by $2j\pi/3$, one obtains
$\gamma_j$. The functions defined by the formulas
\begin{equation}
  \label{Airy}
  w_j(\zeta)=\int_{\gamma_j}e^{-\left(\frac{s^3}3-\zeta
      s\right)}\,ds,\quad j\in\Z_3,\qquad \zeta\in\C,
\end{equation}
are three Airy functions related to the standard
Airy function Ai by the formulas (see, e.g.,~\cite{W})
\begin{equation}
  \label{eq:Ai}
  w_j(\zeta)=2\pi i \omega^j{\rm Ai}(\omega^j \zeta), \quad
  \omega=e^{2\pi i/3}, \qquad \zeta\in\C.
\end{equation}
%
%
%
Assume that $|\arg z| < 2\pi/3$.  As $|z|\to \infty$ one has
\begin{equation}
  \label{ai:as:simple} 
  {\rm Ai}\,(z) =\frac{\exp\left( -\frac23
      z^{3/2}+o(1)\right)} {2\sqrt{\pi}\;z^{1/4}},\
  {\rm Ai}\,(-z) =\frac{\cos\left( \frac23
      z^{3/2}+\frac\pi4+o(1)\right)} {\sqrt{\pi}\;z^{1/4}} \,(1+o(1))
\end{equation}
where we use  the analytic branches of $z\to z^{3/2}$ and $z\to z^{1/4}$ that are
positive for $z>0$ (see~\cite{O:74}, pp. 116, 118 and 392).
\subsection{Notations}
\label{sec:some-agre-notat}
The letter $C$ denotes various positive constants independent
of $z$ and $h$.\\
For two functions $f$ and $g$ defined on a domain $D\subset\C$, we
write that $g(z)=O(f(z))$ in $D$ if $|g(z)|\le C|f(z)|$ for all $z\in D$.\\
\subsection{Solutions in a complex neighborhood of a branch point}
\label{sec:solut-compl-neighb}
\subsubsection{Asymptotic solutions}
\label{sec:asymptotic-solutions-1}
First, let us describe asymptotic solutions
to~\eqref{eq:main}. Therefore, we introduce several objects. 
For a function $f$ defined on $U$, we set
\begin{equation}
  \label{eq:H-def}
  [H(f)](z):=f(z+h)+f(z-h)+v(z)f(z)\quad\text{if}\quad
  \{z-h,z,z+h\}\subset U.
\end{equation}
We let
\begin{equation}\label{g:def}
  g(z):=\frac{\sh\left(\sqrt{\zeta(z)}\zeta'(z)\right)}
  {\sqrt{\zeta(z)}},\qquad z\in U, 
\end{equation}
where the determination of the square roots in the denominator and the
numerator are the same (the definition of $g$ is independent of its
choice). The function $g$ is analytic in $U$. Possibly
reducing $U$ somewhat, we can and do assume that 
\begin{itemize}
\item $g$  does not vanish in $U$.  
\end{itemize}
We further define
\begin{equation}
  \label{eq:A0}
  A_0(z):=\frac1{\sqrt{g(z)}}.
\end{equation}
The function $A_0$ is analytic in $U$.\\
One has
\begin{Th}
  \label{th:as-sol}
  There exist functions $(A_l)_{l\in\N\cup\{0\}}$ and $(B_l)_{l\in\N}$, {\rm (} $A_0$  
  being defined by~\eqref{eq:A0} {\rm )}, 
  all analytic on $U$ and such that, for any $L\in \N\cup\{0\}$ the following holds.
  Let $w$ be one of the Airy functions $w_j$, \ $j\in\Z_3$. If we define
  \begin{equation}
    \label{w-sans-arg}
    w_h(z)=w\left(\zeta(z)/h^{\frac23}\right),\qquad
    w_h'(z)=w'\left(\zeta(z)/h^{\frac23}\right).
  \end{equation}
  and
  \begin{equation}
    \label{as-sol}
    W(z)=h^{\frac13}w_h(z)\;\sum_{l=0}^L
    h^lA_l(z)+h^{\frac23}w_h'(z)\;\sum_{l=1}^L h^lB_l(z),
  \end{equation} 
  then one has
  \begin{equation}
    \label{eq:as-sol}
    H(W)=O\left(h^{L+2+\frac13}w_h\right)+O\left(h^{L+2+\frac23}w_h'\right).
  \end{equation}
\end{Th}
\noindent We call the formal expression
\begin{equation}
  \label{eq:formal-sol}
  h^{\frac13}w_h(z)\;\sum_{l=0}^\infty h^lA_l(z)+h^{\frac23}w_h'(z)\;\sum_{l=1}^\infty h^lB_l(z)
\end{equation}
an {\it asymptotic solution} to~\eqref{eq:main}.

Theorem~\ref{th:as-sol} is proved in section~\ref{s:as-sol}, where we
describe, inter alia, a way to compute the coefficients $(A_l)_l$ and
$(B_l)_l$.

Let us comment on the results of Theorem~\ref{th:as-sol}. First, we note that,
for the differential equation $-\psi''(z)+v(z)\psi(z)=0$, in a
neighborhood of a simple turning point (a point where $v(z)=0$ and
$v'(z)\ne 0$), there are asymptotic solutions of the form~\eqref{eq:formal-sol} 
(with different coefficients $(A_l)_l$,  $(B_l)_l$ and function $\zeta$). 

To justify the Ansatz~\eqref{eq:formal-sol} for the difference equation, 
one has to derive asymptotic formulas of the form
  \begin{equation}
    \label{eq:w-zeta-h}
    \begin{split}
      w_h(z\pm h)=
      f(z)\,w_h(z) \pm  h^{\frac13}\,g(z)\,w_h'(z)+\dots\,,
    \end{split}
  \end{equation}
  where $f(z)=\cosh(\sqrt{\zeta(z)}\zeta'(z))$ and the dots denote
  smaller order terms.  If one tries to prove this formula using
  Taylor expansions for the left hand side, one has to handle an
  infinite number of infinite subsequences of terms of the same
  order. So, an effective resummation of these sequences is
  required. As we see in this paper, to derive formulas analogous 
  to~\eqref{eq:w-zeta-h}, instead of resummation of Taylor series,  
  it is very natural to use tools from complex  analysis.

Formula~\eqref{eq:w-zeta-h} imply that 
\begin{equation}\label{p-and-zeta}
H\left(w_h\right)(z)=\left(2\cosh(\sqrt{\zeta(z)}\zeta'(z))+v(z)\right)w_h(z)+\dots\,.
\end{equation}
In view of~\eqref{eq:momentum} and~\eqref{dif-eq:zeta},
the leading term in the right-hand side of~\eqref{p-and-zeta} is zero.

Finally, we note that  if $h^{-\frac23}|\zeta(z)|$ is large, 
then $w_h(z)$ and $w_h'$ in~\eqref{as-sol} can be  replaced by their 
asymptotic representations. As a result, in view of~\eqref{ai:as:simple}, 
the leading term in~\eqref{as-sol} turns into a linear combination of the leading 
terms from~\eqref{eq:main:wkb}.

\subsubsection{Solutions with standard asymptotic behavior}
\label{sec:exact-solutions-with}
Our main result  is
\begin{Th}
  \label{main:th}
  Let $L\in \N$, and let $W$ be one of the functions constructed in
  Theorem~\ref{th:as-sol} for the order $L$. Then there exists an
  $h$-independent neighborhood $\U\subset U$ of $z_0$ such that, for
  sufficiently small $h$, there exists $\psi$,  a solution
  to equation~\eqref{eq:main} that is analytic in $\U$ and admits there the
  asymptotic representation
  \begin{equation}
    \label{eq:asymptotics-main-w:0}
    \psi(z)=W(z)+O(w_h\,h^{L+1+\frac13})+O(w_h '\,h^{L+1+\frac23})
  \end{equation}
  where $w_h$ and $w_h'$ are defined   in~\eqref{w-sans-arg}. 
\end{Th}
\noindent Theorem~\ref{main:th} is proved in sections~\ref{s:CD}
and~\ref{s:proof:main}.\\
Let us briefly explain the idea of the proof of Theorem~\ref{main:th}.
First, in section~\ref{s:CD}, using the approximate solutions
constructed in Theorem~\ref{th:as-sol}, we construct a parametrix $R$,
i.e., an operator such that, for suitable functions $f$, one has
$HRf=f+Df$, where $H$ is defined in~\eqref{eq:H-def} and $D$ is a
small operator.  The operator $D$ is a singular integral operator. We
estimate its norm using natural geometric objects of the
complex WKB method. This allows us  to prove Theorem~\ref{main:th} on
some special subdomains of $U$. In section~\ref{s:proof:main}, we
study the thus constructed solutions on larger domains and complete
the proof of Theorem~\ref{main:th}.\\
To complete this short description, let us underline that, as
equation~\eqref{eq:main} is non-local in $z$, the ideas of analysis
of~\eqref{eq:main} are different from those used to study the
analogous differential equation.
\subsection{Related results}
\label{rel-res}
The WKB asymptotics of solutions of difference equations on $\Z$ with
``slowly varying'' coefficients have been studied since the end of
1960-s. In~\cite{V-Y:70} and~\cite{T:74}, the authors essentially
studied equations of the form
\begin{equation}
  \label{eq-on-Z}
  Y_{k+1}=M(hk) Y_k,\quad k\in \Z,
\end{equation}
with a small positive $h$ and an $(n\times n)$-matrix valued function $M$
defined on $\R$.  We note that if
\begin{equation}
  \label{eq-on-R}
  Y(x+h)=M(x) Y(x),\quad x\in \R,
\end{equation}
then, the sequence $(Y_k)_{k\in\Z}=(Y(kh))_{k\in\Z}$
satisfies~\eqref{eq-on-Z}.  We note also that equation~\eqref{eq:main}
restricted to $\R$ is equivalent to~\eqref{eq-on-R} with
$M(x)=\begin{pmatrix} -v(x) &-1 \\ 1& 0\end{pmatrix}$, and that   a
turning point  for equation~\eqref{eq:main}
is a point $x$ where the eigenvalues of the matrix $M(x)$ coincide. \\
The short note~\cite{V-Y:70} is essentially devoted to the case where
all the eigenvalues of the matrix $M$ in~\eqref{eq-on-Z} are distinct.
In~\cite{T:74} the author constructed asymptotic solutions
to~\eqref{eq-on-Z} in a small (depending on $h$) neighborhood of a
point where two eigenvalues of $M(x)$ become equal.\\
In~\cite{CC:96} the authors considered difference equations
of the form 
\begin{equation*}
   \sum_{j=I}^Ja_j(hk,h)\,y_{k+j}=0, \quad k\in\Z.
\end{equation*}
We note that this class includes the difference Schr{\"o}dinger equations
\begin{equation*}
  y_{k+1}+y_{k-1}+v(hk) y_k=0, \quad k\in\Z.
\end{equation*}
The authors described the asymptotics of solutions to this equations
for $hk$ being in a small (as $h\to 0$) neighborhood of a point where
$v(x)\in\{\pm 1\}$.\\
We mention also three (series of) papers motivated  by problems originating in 
the theory of orthogonal polynomials.\\
First, there is  a series of papers by J.S. Geronimo and co-authors,
see, e.g.~\cite{G-at-al} and references therein, devoted to uniform
asymptotic formulas for solutions to the equation
$a_{k+1} \psi_{k+1} + b_k \psi_k + a_k \psi_{k-1} = \lambda\psi_k$, \
$k\in \Z$, where $\lambda$ is the spectral parameter, and the coefficients
$a_k$ are positive and $b_k$ are real.  \\
Also we mention   papers by R.Wong and coauthors, see e.g.,~\cite{W:03}, 
who also studied solutions to three terms recurrence relations  with real 
coefficients for large values of the integer variable.\\
Finally, we mention paper~\cite{Dobro} where the authors
studied WKB asymptotics of solutions to a difference equation 
using the Maslov canonical operator.\\
There are more papers devoted to the subject. The reader can find more references
in the papers that we mentioned above. 
 
To the best of our knowledge, the present paper is the first where one
rigorously obtains uniform asymptotics of analytic solutions to a
difference equation on $\C$ in an $h$-independent neighborhood of a
turning point.
\section{The space of solutions to equation~\eqref{eq:main}}
\label{s:space-of-sol}
The observations that we discuss now are well-known in the theory of
difference equations and are easily proved.\\
Let $c\in \R$ and $I=\{z\in U: \im z= c\}$. We assume that the length of the
segment $I$ is sufficiently large (with respect to $h$) and discuss
the set $S$ of solutions to equation~\eqref{eq:main} on $I$.\\
Let $\{f,g\}\subset S$. The expression
\begin{equation}
  \label{eq:wronskian}
  (f(z),g(z))=f(z+h)g(z)-f(z)g(z+h), \quad z,\,z+h\in I,
\end{equation}
is called {\it the Wronskian} of $f$ and $g$.   It is $h$-periodic in $z$.\\
If the Wronskian of two solutions does not vanish, they form a basis
in $S$, i.e, $\psi\in S$ if and only if
\begin{equation}\label{eq:three-solutions}
  \psi(z)=a(z)f(z)+b(z)g(z),\quad z,\,z+h\in I,
\end{equation}
where $a$ and $b$ are $h$-periodic complex coefficients. One has
\begin{equation}
  \label{eq:periodic-coef}
  a(z)=\frac{(\psi(z),\,g(z))}{(f(z),\,g(z))},\quad 
  b(z)=\frac{(f(z),\,\psi\,(z))}{(f(z),\,g(z))}.
\end{equation} 
\section{Asymptotic solutions: the proof of 
Theorem~\ref{th:as-sol}}
\label{s:as-sol}
\subsection{The proof of Theorem~\ref{th:as-sol} up to two
  propositions}
\label{sec:proof-theor-refth}
First, we formulate two statements needed to construct asymptotic
solutions to~\eqref{eq:main}. Below we use the notations introduced
in~\eqref{w-sans-arg}.
\begin{Pro}
  \label{pro:Ahw}
  Let $A$ be analytic in $U$. Let $N\in\N$. If $\{z-h,z, z+h\}\subset  U$,
  \begin{equation}
    \label{eq:Hv}
    \begin{split}
      H\left(A\,h^{\frac13}w_h\,\right)&= h^{\frac13}w_h\;\sum_{l=2}^N
      h^la_l+O\left(h^{N+1+\frac13}w_h\right)\\
      &\hspace{3cm}+h^{\frac23}w_h'\;\sum_{l=1}^N
      h^lb_l+O\left(h^{N+1+\frac23}w_h'\right)
    \end{split}
  \end{equation}
   as $h\to 0$.   All the coefficients $(a_l)_{l\ge 2}$ and
  $(b_l)_{l\ge 1}$ are analytic in $U$, independent of the
  choice of $w$ in~\eqref{w-sans-arg}, and
  \begin{equation}
    \label{b1}
    b_1=b_1[A]=Ag\frac{d}{dz}\log\left(A^2g\right).
  \end{equation}
\end{Pro}
\begin{Pro}
  \label{pro:Bhw}
  Let $B$ be analytic in $U$, and let  $N\in\N$.
  If $\{z-h,z, z+h\}\subset  U$,   
  \begin{equation}
    \label{eq:Hv-prime}
    \begin{split}
      H\left(B\,h^{\frac23}w_h'\right)&= h^{\frac13}w_h\;
      \sum_{l=1}^N h^lc_l + O\left(h^{N+1+\frac13}w_h\right)\\
      &\hspace{3cm} +h^{\frac23}w_h'\,\sum_{l=2}^N
      h^ld_l+O\left(h^{N+1+\frac23}w_h'\right)
    \end{split}
  \end{equation}
  as $h\to 0$. All the coefficients $(c_l)_{l\ge 1}$ and
  $(d_l)_{l\ge 2}$ are analytic in $U$, independent of the
  choice of $w$ in~\eqref{w-sans-arg}, and
  \begin{equation}
    \label{c1}
    c_1=c_1[B]=\zeta Bg\frac{d}{dz}\log(\zeta B^2g).
  \end{equation}
\end{Pro}
\noindent Before proving Propositions~\ref{pro:Ahw} and~\ref{pro:Bhw},
we use them to prove Theorem~\ref{th:as-sol}. The proof is done by
induction on the order $L$. For $L=0$, one has $W=A_0h^{\frac13}w_h$.
In view of~\eqref{eq:A0} and~\eqref{b1}, the coefficient $b_1$
corresponding to $A=A_0$ is equal to $0$. So, the statement of
Theorem~\ref{th:as-sol} for $L=0$ immediately follows
from~\eqref{eq:Hv} with $N=1$.\\
Now, we assume that Theorem~\ref{th:as-sol} is proved up to the order
$L=L_0-1$, $L_0\in\N$. Let us prove it for $L=L_0$. We set
\begin{equation}
  \label{as-sol:1}
  W(z)=h^{\frac13}w_h(z)\;\sum_{l=0}^{L_0} h^lA_l(z)
  +h^{\frac23}w_h'(z)\;\sum_{l=1}^{L_0} h^lB_l(z),
\end{equation}
where $(A_l)_{l<L_0}$ and $(B_l)_{l<L_0}$ are chosen as in the case
$L=L_0-1$, $A_{L_0}$ and $B_{L_0}$ still having to be chosen.  By the
induction hypothesis
\begin{multline}
  \label{eq:as-sol:1}
  H\left(h^{\frac13}w_h\;\sum_{l=0}^{L_0-1} h^lA_l
    +h^{\frac23}w_h'\;\sum_{l=1}^{L_0-1} h^lB_l\right)\\=
  O\left(h^{L_0+1+\frac13}w_h\right)+O\left(h^{L_0+1+\frac23}w_h'\right).
\end{multline}
In view of Propositions~\ref{pro:Bhw} and~\ref{pro:Ahw} this implies
that
\begin{equation*}
  \label{eq:as-sol:2}
  \begin{split}
    H\left(h^{\frac13}w_h\sum_{l=0}^{L_0-1} h^lA_l
      +h^{\frac23}w_h'\sum_{l=1}^{L_0-1} h^lB_l\right)=&
    ah^{L_0+1+\frac13}w_h+bh^{L_0+1+\frac23}w_h'\\
    &+O\left(h^{L_0+2+\frac13}w_h\right)+O\left(h^{L_0+2+\frac23}w_h'\right),
  \end{split}
\end{equation*}
where $a$ and $b$ are analytic functions in $U$.  On the other hand,
using~\eqref{eq:Hv} and~\eqref{eq:Hv-prime} with $N=1$, we get
\begin{gather*}
  H(A_{L_0}h^{L_0+\frac13}w_h)=h^{L_0+1+\frac23}w_h'\,b_1[A_{L_0}]+
  O(h^{L_0+2+\frac13}w_h)+O(h^{L_0+2+\frac23}w_h'),\\
  H(B_{L_0}h^{L_0+\frac23}w_h')=h^{L_0+1+\frac13}w_h\,c_1[B_{L_0}]+
  O(h^{L_0+2+\frac13}w_h)+O(h^{L_0+2+\frac23}w_h').
\end{gather*}
Therefore,
\begin{equation*}
  \begin{split}
    HW&=h^{L_0+1}\left( h^{\frac13}w_h\,(a+c_1[B_{L_0}])+
      h^{\frac23}w_h'\,(b+b_1[A_{L_0}])\right)\\
    &\hspace{4cm} +O(h^{L_0+2+\frac13}w_h)+O(h^{L_0+2+\frac23}w_h').
  \end{split}
\end{equation*}
So, to prove Theorem~\ref{th:as-sol}, it suffices to choose $A_{L_0}$
and $B_{L_0}$ so that
\begin{equation*}
  a+c_1[B_{L_0}]=0,\quad\text{and}\quad  b+b_1[A_{L_0}]=0.
\end{equation*}
In view of~\eqref{b1} and~\eqref{c1}, these relations are equivalent
to the equations
\begin{equation}\label{eq:AL0,BL0}
  \zeta B_{L_0}g\frac{d}{dz}\log\left(\zeta B_{L_0}^2g\right)=-a,
  \quad\text{and}\quad
  A_{L_0}g\frac{d}{dz}\log\left( A_{L_0}^2g\right)=-b.
\end{equation}
One constructs solutions to these equations by the formulas
\begin{equation}\label{AL0,BL0}
  A_{L_0}(z)=-\frac1{2\sqrt{g(z)}}\int_{z_0}^z\frac{b\,dz}{\sqrt{g}},
  \quad\text{and}\quad
  B_{L_0}(z)=-\frac1{2\sqrt{\zeta(z)g(z)}}\int_{z_0}^z\frac{a\,dz}{\sqrt{\zeta g}}.
\end{equation}
As
\begin{itemize}
\item $g$ and $\zeta$ are analytic in $U$,
\item $g$ does not vanish in $U$,
\item $\zeta$ vanishes in $U$ only at $z_0$ where it has a simple
  zero,
\end{itemize}
the coefficients $A_{L_0}$ and $B_{L_0}$ are analytic in $U$.  This
completes the proof of Theorem~\ref{th:as-sol}.  \qed
\begin{Rem}
  The function $B_{L_0}$ constructed by~\eqref{AL0,BL0} is the only
  solution to the first equation in~\eqref{eq:AL0,BL0} that is
  analytic in $U$. The function $A_{L_0}$ constructed
  by~\eqref{AL0,BL0} is unique up to a solution to the homogeneous
  equation $A_{L_0}g\frac{d}{dz}\log\left( A_{L_0}^2g\right)=0$ that
  is proportional to $A_0$ given by~\eqref{eq:A0}.
\end{Rem}
\subsection{The proof of Proposition~\ref{pro:Ahw}}
\label{sec:proof-prop-refpr}
Consider $(w_j)_{j\in\Z_3}$ the three  solutions
to the Airy equation~\eqref{eq:Airy} defined by~\eqref{Airy}.  Let $w$
be $w_j$ for some $j\in\Z_3$ and let $\gamma$ be the corresponding
integration path $\gamma_j$ in~\eqref{Airy}.\\
Note that
\begin{equation}
\label{eq:h-airy}
h^{\frac13}w(h^{-\frac23}\zeta)=\int_\gamma e^{-\frac1h\left(\frac{t^3}3-t\zeta\right)}\,dt,
\qquad
h^{\frac23}w'(h^{-\frac23}\zeta)=\int_\gamma e^{-\frac1h\left(\frac{t^3}3-t\zeta\right)}t\,dt.
\end{equation}
Below, we use the notations from~\eqref{w-sans-arg}.  Let $K\subset U$ is a
closed disk centered at $z_0$ and independent of $h$.  Below, we
assume that $z\in K$ and $h$ is sufficiently small.  The proof of the
asymptotics of $H(Ah^{\frac13}w)$ in $K$ as $h\to 0$ is broken into
several steps.
\\
{\bf 1.} \ In view of~(\ref{eq:h-airy}), we get
\begin{gather}
  \label{Aw:int}
  H\left(Ah^{\frac13}w\right)=
  \int_{\gamma}e^{-\frac1h\left(\frac{t^3}3-t\zeta(z)\right)}F_0(t,z,h)\,dt,\\
  \label{Aw:F}
  F_0(t,z,h)=A(z+h)e^{\frac{t}h(\zeta(z+h)-\zeta(z))}+
  A(z-h)e^{\frac{t}h(\zeta(z-h)-\zeta(z))}+v(z)A(z).
\end{gather}
Note that $(t,z,h)\mapsto F_0(t,z,h)$ is analytic in $\C\times K\times V$, where
$V$ is a sufficiently small neighborhood of zero.
\\
{\bf 2.} \ To get the asymptotics of the integral in~\eqref{Aw:int},
we apply the well-known method described in detail in section 4 of
chapter VII of~\cite{W}. First, we represent $F_0(t,z,h)$ in the form
\begin{equation}
  \label{F-f0}
  F_0(t,z,h)=a_0(z,h)+b_0(z,h)t+(t^2-\zeta(z)) f_0(t,z,h)
\end{equation}
with
\begin{gather}
  \label{eq:a0}
  a_0(z,h)=\frac12\,\left(F_0(\sqrt{\zeta(z)},z,h)+
    F_0(-\sqrt{\zeta(z)},z,h)\right),
  \\
  \label{eq:b0}
  b_0(z,h)=\frac1{2\sqrt{\zeta(z)}}\,
  \left(F_0(\sqrt{\zeta(z)},z,h)-F_0(-\sqrt{\zeta(z)},z,h)\right),
\end{gather}
where, in~\eqref{eq:a0} and~\eqref{eq:b0}, we use one and the same branch of
$\sqrt{\zeta(z)}$.
\\
Both $a_0$ and $b_0$ are analytic in $(z,h)\in K\times V$ (we remove
the removable singularities at $z=0$).  With $a_0$ and $b_0$ so chosen,
it is easily seen that the function $f_0$ is analytic 
in $(t,z,h)\in \C\times K\times V$.
\\
{\bf 3.} \ Substituting~\eqref{F-f0} into~\eqref{Aw:int} and
integrating by parts, we get
\begin{equation}
  \label{int-of-F1}
  \begin{split}
    H\left(Ah^{\frac13}w_h\right)= a_0h^{\frac13}w_h
    +b_0h^{\frac23}w_h' +\int_\gamma
    e^{-\frac1h\left(\frac{t^3}3-t\zeta(z)\right)}(t^2-\zeta(z))
    f_0\,dt\\
    =a_0h^{\frac13}w_h +b_0h^{\frac23}w_h' +h\int_\gamma
    e^{-\frac1h\left(\frac{t^3}3-t\zeta(z)\right)} F_1(t,z,h)\,dt.
  \end{split}
\end{equation}
where $F_1(t,z,h)=\frac{\partial f_0}{\partial t}(t,z,h)$.
\\
{\bf 4.} \ Now, we transform the last integral in~\eqref{int-of-F1},
the one containing $F_1$, in the same way as we transformed the
integral with $F_0$ from~\eqref{Aw:int}.
\\
For a fixed positive integer $N$, we repeat this procedure inductively $N+2$
times. Reasoning as above, one proves that
\begin{gather}
  \label{pre-HAw}
  H\left(A\,h^{\frac13}w_h\,\right)= h^{\frac13}w_h\;\sum_{l=0}^{N+1}
  h^l a_l(z,h)+
  h^{\frac23}w_h'\;\sum_{l=0}^{N+1} h^l b_l(z,h)+h^{N+2}I_{N+2},\\
  \intertext{where, for $l\in\N\cup\{0\}$, we have defined}
  \label{Il-def}
  I_{l}=\int_\gamma e^{-\frac1h\left(\frac{t^3}3-t\zeta(z)\right)}
  F_l(t,z,h)\,dt.
\end{gather}
As when $l=0$, the coefficients $a_{l}$ and $b_{l}$ are expressed in
terms of $F_l$ by
\begin{gather*}
  a_l(z,h)=\frac12\,\left(F_l(\sqrt{\zeta(z)},z,h)+
    F_l(-\sqrt{\zeta(z)},z,h)\right),
  \\
  b_l(z,h)=\frac1{2\sqrt{\zeta(z)}}\,
  \left(F_l(\sqrt{\zeta(z)},z,h)-F_l(-\sqrt{\zeta(z)},z,h)\right),
\end{gather*}
and the function $f_l$ is defined by the relation
\begin{equation}
  \label{eq:3}
  F_l(t,z,h)=a_l(z,h)+b_l(z,h)t+(t^2-\zeta(z)) f_l(t,z,h)
\end{equation}
Finally, for $l\geq1$, one has
$F_l(t,z,h)=\frac{\partial f_{l-1}}{\partial t}(t,z,h)$.\\
For $l\in\N\cup\{0\}$, the coefficients
$a_l$ and $b_l$ are analytic in $(z,h)\in K\times V$.\\
{\bf 5.} \ To estimate the integrals $I_l$, one has to estimate the
functions $F_l$. Below, the constants $C$ are  independent on $z$,
$h$ and $t$. The symbol $O(\cdot)$ is subsequently used for estimates
uniform in $z$, $t$ and $h$.\\

Let us assume that  that$(t,z,h)\in \C\times K\times V$ and show that there 
exists a constant $C_0>0$ such that, for any $l\in\N$, one has
\begin{equation}
  \label{eq:Fl-est}
  F_l(t,z,h)=O\left(e^{C_0|t|}\right), 
\end{equation}
where the implicit constant in~\eqref{eq:Fl-est} depends only on the
index $l$. \\
This estimate is obvious for $F_0$. Let us assume that it is proved
for some $l=l_0$ and prove it for $l=l_0+1$.\\
Clearly, $\zeta(z)$ is bounded on $K$.  In view of the definitions of
$(a_l)_{l\ge 0}$, $(b_l)_{l\ge 0}$ and the induction
hypothesis, we have $a_{l_0}(z,h)=O(1)$ and $b_{l_0}(z,h)=O(1)$.
These observations, the definition of $f_l$~\eqref{eq:3} and the
induction hypothesis imply that there exists $R>0$ independent of $h$
such that, for all $|t|\ge R$, $f_{l_0}(t,z,h)= O(e^{C_0|t|})$. By the
maximum principle, this implies that $f_{l_0}$ satisfies this estimate
for all $t\in \C$.  Now the Cauchy estimates for the derivatives of the
analytic functions imply~\eqref{eq:Fl-est} for $l=l_0+1$.
\\
{\bf 6.} \ Let us prove that
\begin{equation}
  \label{eq:Il-est}
  I_l=O\left(h^{\frac13}w_h\right)+O\left(h^{\frac13}w_h'\right).
\end{equation}
Therefore, one essentially has to repeat the reasoning made in Section
4, Chapter VII of~\cite{W}. So, we omit some details. \\
If $Z=h^{-\frac23}\zeta(z)$ is bounded by a constant, setting
$T=h^{-\frac13}t$, we change variable in~\eqref{Il-def}. In view of
step 5, we get
\begin{equation*}
  I_l=h^{\frac13}\int_\gamma e^{-\left(\frac{T^3}3-TZ\right)}
  O(e^{C_0h^{\frac13}|T|})\,dT =O(h^{\frac13}),
\end{equation*}
and this leads to~\eqref{eq:Il-est} as $w$ and $w'$ have no common
zero.\\
If $Z=h^{-\frac23}\zeta(z)$ is large, we estimate the integral $I_l$
using the method of steepest descent. In view of the fifth step, we
have
\begin{equation*}
  I_l=\int_\gamma e^{-\frac1h\left(\frac{t^3}3-t\zeta(z)\right)}
  O(e^{C_0|t|})\,dt. 
\end{equation*}
We deform the integration path to a path of steepest descent for
$e^{-\frac1h\left(\frac{t^3}3-t\zeta(z)\right)}$ exactly as when
computing the asymptotics of the Airy function $w$, i.e., the
asymptotics of the integral
$\int_\gamma e^{-\frac1h\left(\frac{t^3}3-t\zeta(z)\right)} \,dt$. The
saddle points $\pm \sqrt{\zeta(z)}$ are uniformly bounded when
$z\in K$. Let $r>0$ be sufficiently large for the saddle points to be
inside the disk of radius $r$ centered at $0$. We compute the
asymptotics of the integral over $\gamma\cap \{|t|\le r\}$ directly by
means of the method of steepest descents and, comparing the answer
with the asymptotics of the Airy function $w(Z)$ as $Z\to\infty$, we find that this
integral is bounded by $O(h^{\frac13}w_h)+O(h^{\frac23}w_h')$. The
integral over the remaining part of $\gamma$ quickly tends to $0$ as
$h\to0$: actually, it is exponentially small with respect to
$O(h^{\frac13}w_h)+O(h^{\frac23}w_h')$. This yields~\eqref{eq:Il-est}.
\\
{\bf 7.} \ Formula~\eqref{pre-HAw} and estimate~\eqref{eq:Il-est} lead
to the representation
\begin{equation}
 \label{HAw:before-Taylor}
  \begin{split}
    H\left(A\,h^{\frac13}w_h\,\right)(z)&= h^{\frac13}w_h(z)\;
    \sum_{l=0}^{N} h^l a_l(z,h)+
    h^{\frac23}w_h'(z)\;\sum_{l=0}^{N} h^l b_l(z,h)\\
    &+O\left(h^{N+1+\frac13}w_h(z)\right)+O\left(h^{N+1+\frac23}w_h'(z)\right),
  \end{split}
\end{equation}
The coefficients $(a_l)_{l\in\N\cup\{0\}}$ and
$(b_l)_{l\in\N\cup\{0\}}$ being analytic in $h$, we can approximate
them by Taylor polynomials. This yields
\begin{equation}
  \label{almost-the-result}
  \begin{split}
    H\left(A\,h^{\frac13}w_h\,\right)(z)&= h^{\frac13}w_h(z)\;\sum_{l=0}^{N}
    h^l a_l(z)+
    h^{\frac23}w_h'(z)\;\sum_{l=0}^{N} h^l b_l(z)\\
    &+O\left(h^{N+1+\frac13}w_h(z)\right)+O\left(h^{N+1+\frac23}w_h'(z)\right),
  \end{split}
\end{equation}
where $(a_l(z))_{l\in\N\cup\{0\}}$ and $(b_l(z)))_{l\in\N\cup\{0\}}$ are new
coefficients independent of $h$. In particular, one has
\begin{equation}\label{a-et-b}
  \begin{split}
    a_0(z)=a_0(z,0),\qquad& b_0(z)=b_0(z,0),\\
    a_1(z)=a_1(z,0)+\frac{\partial a_0}{\partial h}(z,0),\qquad
    &b_1(z)=b_1(z,0)+\frac{\partial b_0}{\partial h}(z,0).
  \end{split}
\end{equation}
Now, to complete the proof of Proposition~\ref{pro:Ahw}, it suffices
to compute $a_0$, $b_0$, $a_1$, $b_1$.\\
{\bf 8.} \ Let us check that
\begin{equation}\label{a0:Taylor}
  a_0(z,0)=\frac{\partial a_0}{\partial h}(z,0)=0.
\end{equation}
Substituting~(\ref{Aw:F}) into~(\ref{eq:a0}), we get
\begin{equation*}
  \begin{split}
    a_0(z,h)=A(z+h)\ch&\left(\sqrt{\zeta(z)}\;\frac{\zeta(z+h)-\zeta(z)}h\right)\\
    &
    +A(z-h)\ch\left(\sqrt{\zeta(z)}\;\frac{\zeta(z-h)-\zeta(z)}h\right)+v(z)A(z).
  \end{split}
\end{equation*}
Thus,
\begin{equation*}
  a_0(z,0)=A(z)\left(2\,\ch(\sqrt{\zeta(z)}\zeta'(z))+v(z)\right).
\end{equation*}
Recall that the complex momentum $p$ is defined
in~(\ref{eq:momentum}). In view of~(\ref{dif-eq:zeta}) we get
\begin{equation}
  \label{eq:zeta-v}
  2\ch\left(\sqrt{\zeta(z)}\zeta'(z)\right)+v(z)=0.
\end{equation}
So, $a_0(z,0)=0$. As $a_0(z,h)$ is even in $h$, we also see that
$\frac{\partial a_0}{\partial h}(z,0)=0$.
\\
{\bf 9.} \ Let us check that
\begin{equation}\label{b0:Taylor}\dsize
  b_0(z,0)=0,\quad \frac{\partial b_0}{\partial h}(z,0)=
  2A'(z)\frac{\sh\left(\sqrt{\zeta(z)}\zeta'(z)\right)}{\sqrt{\zeta(z)}}+
  A(z)\ch\left(\sqrt{\zeta(z)}\zeta'(z)\right)\zeta''(z).
\end{equation}
Substituting~(\ref{Aw:F}) into~(\ref{eq:b0}), we get
\begin{equation*}\dsize
  b_0(z,h)=A(z+h)\frac{\sh\left(\sqrt{\zeta(z)}\;\frac{\zeta(z+h)-\zeta(z)}h\right)}
  {\sqrt{\zeta(z)}}
  +A(z-h)\frac{\sh\left(\sqrt{\zeta(z)}\;\frac{\zeta(z-h)-\zeta(z)}h\right)}
  {\sqrt{\zeta(z)}}.
\end{equation*}
Clearly, $b_0(z,h)$ is odd in $h$, and so $b_0(z,0)=0$.  Computing
$\frac{\partial b_0}{\partial h}(z,0)$, we complete the proof
of~\eqref{b0:Taylor}.\\
{\bf 10.} \ To compute $a_1(z,0)$ and $b_1(z,0)$, we first study
$f_0$. Let $r>0$ be such that $|\zeta(z)|\le r^2/2$ for all $z\in K$.  Let
$|t|=r$, $z\in K$ and $h\in V$. Formulas~(\ref{Aw:F}) and~(\ref{eq:zeta-v}) imply
that
\begin{equation*}\dsize
  F_0(t,z,h)=F_0(t,z)+O(h), \ \  
  F_0(t,z)=2A(z)\left(\ch\left(t\zeta'(z)\right)-
    \ch(\sqrt{\zeta(z)}\zeta'(z))\right).
\end{equation*}
This result, the formulas $a_0(z,0)=b_0(z,0)=0$ (see steps 9--10)
and~(\ref{F-f0}) imply that for $|t|=r$  one has
\begin{equation*}
  f_0(t,z,h)=\frac{F_0(t,z)}{t^2-\zeta(z)} +O(h).
\end{equation*}
By the maximum principle for analytic functions, this representation
remains true for all $|t|\le r$.
\\
{\bf 11.} \ The result of the previous step and the Cauchy estimates
for the derivatives of the analytic functions imply that, for
$|t|\le r/2$, $z\in K$ and $h\in V$,  one has
\begin{equation}
  F_1(t,z,h)=F_{1}(t,z)+O(h),\quad\text{where}\quad  F_{1}(t,z)=
  \frac{\partial}{\partial t}\left(\frac{F_{0}(t,z)}{t^2-\zeta(z)}\right).
\end{equation}
Therefore
\begin{equation*}\textstyle
  a_1(z,0)=\frac{F_1\left(\sqrt{\zeta(z)},z\right)+
    F_1\left(-\sqrt{\zeta(z)},z\right)}2, \quad 
  b_1(z,0)=\frac{F_1\left(\sqrt{\zeta(z)},z\right)-
    F_1\left(-\sqrt{\zeta(z)},z\right)}
  {2\sqrt{\zeta(z)}}.
\end{equation*}
As $t\mapsto F_1(t,z)$ is odd, one has
\begin{equation*}
  a_1(z,0)=0,\qquad b_1(z,0)=\frac1t\left.
    \frac{\partial}{\partial t}
    \left(\frac{F_{0}(t, z)}{t^2-\zeta(z)}\right)\right|_{t=\sqrt{\zeta}}.
\end{equation*}
Elementary calculations yield
\begin{equation*}
  b_1(z,0)=\frac{A\zeta'}{2\zeta}\left(\zeta'(z)
    \ch\left(\sqrt{\zeta}\zeta'\right)-
    \frac{\sh\left(\sqrt{\zeta}\zeta'\right)}{\sqrt{\zeta}}\right),\quad 
  \zeta=\zeta(z). 
\end{equation*}

The results of the steps 8, 9 and 11 imply that
\begin{equation*}
  a_0(z)=b_0(z)=a_1(z)=0,\qquad
  b_1(z)=A(z) g(z)\frac{d\log(A^2g)}{dz}(z),
\end{equation*}
where $g$ is the defined in~(\ref{g:def}). Substituting these formulae
into~\eqref{almost-the-result}, we obtain~(\ref{eq:Hv})
and~(\ref{b1}).  This completes the proof of
Proposition~\ref{pro:Ahw}.  \qed
\subsection{Proof of Proposition~\ref{pro:Bhw}}
\label{sec:proof-prop-refpr-1}
The proof of Proposition~\ref{pro:Bhw} being parallel to that of
Proposition~\ref{pro:Ahw}, we concentrate only on the differences and
omit details.\\
We assume that $B$ is analytic in $U$. Let $w$, $\gamma$ and $K$
be as in the proof of Proposition~\ref{pro:Ahw}. We use the notations
from~\eqref{w-sans-arg}. We assume that $z\in K$ and that $h$ is
sufficiently small.  The derivation of the asymptotics of
$H(Bh^{\frac23}w')$ is split into several steps.
\\
{\bf 1.} \ We get
\begin{gather}\label{Bw:int}
  H\left(Bh^{\frac23}w_h'\right)=
  \int_{\gamma}e^{-\frac1h\left(\frac{t^3}3-t\zeta(z)\right)}
  G_0(t,z,h)\,dt,\\
  \label{Bw:G}
  G_0(t,z,h)=t\left(B(z+h)e^{\frac{\zeta(z+h)-\zeta(z)}ht}+
    B(z-h)e^{\frac{\zeta(z-h)-\zeta(z)}ht}+v(z)B(z)\right).
\end{gather}
\\
{\bf 2.} \ Fix an $N\in \N$. Reasoning as in steps 1-6 of the proof of
Proposition~\ref{pro:Ahw}, instead of~\eqref{HAw:before-Taylor}, for
$z\in K$ and for sufficiently small $h$, we prove that
\begin{equation*}
  \label{HBw:before-Taylor}
  \begin{split}
    H\left(B\,h^{\frac23}w_h'\right)= h^{\frac13}w_h\;\sum_{l=0}^{N}
    h^l c_l(z,h)&+
    h^{\frac23}w_h'\;\sum_{l=0}^{N} h^l d_l(z,h)\\
    &+O\left(h^{N+1+\frac13}w_h\right)+O\left(h^{N+1+\frac23}w_h'\right),
  \end{split}
\end{equation*}
where, for $l\in\N\cup\{0\}$, one computes
\begin{gather}
  \label{eq:cl}
  \textstyle c_l(z,h)=\frac12\,\left(G_l(\sqrt{\zeta(z)},z,h)
    +G_l(-\sqrt{\zeta(z)},z,h)\right),
  \\
  \label{eq:dl}
  d_l(z,h)=\frac1{2\sqrt{\zeta(z)}}\,
  \left(G_l(\sqrt{\zeta(z)},z,h)-G_l(-\sqrt{\zeta(z)},z,h)\right),
  \\
  \label{gl-GL}
  g_l(t,z,h)=\frac{G_l(t,z,h)-c_l-d_lt}{t^2-\zeta(z)},\quad
  G_{l+1}=\frac{\partial g_l}{\partial t}.
\end{gather}
In~\eqref{eq:cl} and~\eqref{eq:dl}, we use one and the same 
branch of $\sqrt{\zeta(z)}$.\\
Approximating the $(c_l)_{l\in\N\cup\{0\}}$ and
$(d_l)_{l\in\N\cup\{0\}}$ as functions of $h$ by Taylor polynomials,
we get
\begin{equation}
  \label{B:almost-the-result}
  \begin{split}
    H\left(B\,h^{\frac23}w_h'\right)= h^{\frac13}w_h\;\sum_{l=0}^{N}
    h^l c_l(z)&+
    h^{\frac23}w_h'\;\sum_{l=0}^{N} h^l d_l(z)\\
    &+O\left(h^{N+1+\frac13}w_h\right)+O\left(h^{N+1+\frac23}w_h'\right).
  \end{split}
\end{equation}
One has
\begin{equation}
  \label{c-et-d}
  \begin{split}
    c_0(z)=c(z,0),\qquad& d_0(z)=d(z,0),\\
    c_1(z)=c_1(z,0)+\frac{\partial c_0}{\partial h}(z,0),\qquad
    &d_1(z)=d_1(z,0)+\frac{\partial d_0}{\partial h}(z,0).
  \end{split}
\end{equation}
\\
{\bf 3.} \ Substituting~(\ref{Bw:G}) into formula~(\ref{eq:dl}) with
$l=0$, we obtain the formulas
\begin{equation}\label{d0:Taylor}
  d_0(z,0)=\frac{\partial d_0}{\partial h}(z,0)=0
\end{equation}
in the same way as we obtained~\eqref{a0:Taylor}.
\\
{\bf 4.} \ Let us check that
\begin{equation}\label{c0:Taylor}\textstyle
  c_0(z,0)=0,\quad \frac{\partial c_0}{\partial h}(z,0)=\sqrt{\zeta}
  \left(2B'(z)\sh\left(\sqrt{\zeta}\zeta'\right)+
    B(z)\ch\left(\sqrt{\zeta}\zeta'\right)\sqrt{\zeta}\zeta''\right),
\end{equation}
where $\zeta=\zeta(z)$. Substituting~(\ref{Bw:G}) into
formula~(\ref{eq:cl}) with $l=0$, we get
\begin{equation*}
  \begin{split}
    c_0(z,h)=B(z+h)&\sqrt{\zeta(z)}\;
    \sh\left(\sqrt{\zeta(z)}\;\frac{\zeta(z+h)-\zeta(z)}h\right)\\
    &+B(z-h)\sqrt{\zeta(z)}\;
    \sh\left(\sqrt{\zeta(z)}\;\frac{\zeta(z-h)-\zeta(z)}h\right).
  \end{split}
\end{equation*}
The coefficient $c_0(z,h)$ is odd in $h$, and so $c_0(z,0)=0$.
Computing $\frac{\partial c_0}{\partial h}(z,0)$, we complete the
proof of~\eqref{c0:Taylor}.
\\
{\bf 5.} \ As $G_1$ is computed in the same way as $F_1$ in steps
11--12 of the proof of Proposition~\ref{pro:Ahw}, we omit the details
and write down the result:
\begin{equation}\label{G1:as}
  G_1(t,z,h)=G_{1}(t,z)+O(h),\quad  G_{1}(t,z)=2B
  \frac{\partial}{\partial t}\frac
  {t\left(\ch\left(t\zeta'\right)-\ch\left(\sqrt{\zeta}\zeta'\right)\right)}
  {t^2-\zeta}, 
\end{equation}
where $\zeta=\zeta(z)$. This representation is locally uniform in $t$
and uniform in $z\in K$.
\\
{\bf 6.} Having described $G_1$, we easily get the formulas
\begin{equation*}
  d_1(z,0)=0,\quad c_1(z,0)=\frac{B \zeta'}2
  \left(\zeta'\ch\left(\sqrt{\zeta}\zeta'\right)+
    \frac{\sh\left(\sqrt{\zeta}\zeta'\right)}{\sqrt{\zeta}}\right),\quad 
  \zeta=\zeta(z).
\end{equation*}
We again omit elementary details and only note that the first formula
follows from the evenness of the function $t\mapsto G_1(t,z)$.
\\
{\bf 7.} \ Substituting the results of steps 3,4 and 6
into~\eqref{c-et-d}, we get
\begin{equation*}
  c_0(z)=d_0(z)=d_1(z)=0,\qquad
  c_1(z)=\zeta(z) B(z) g(z)\frac{d\log(\zeta B^2g)}{dz}(z),
\end{equation*}
where $g$ is the function from~(\ref{g:def}). Substituting these
formulas into~\eqref{B:almost-the-result}, we prove the statement of
Proposition~\ref{pro:Bhw}.  \qed
\section{Properties of asymptotic solutions}
\label{sec:prop-asympt-solut}
We now study basic properties of the asymptotic solutions.  More
precisely, we fix an integer $L$ and study the functions
$(W_j)_{j\in\Z_3}$, i.e., the functions $W$ from
Theorem~\ref{th:as-sol} corresponding to the chosen $L$ and to the
Airy functions $w=w_j$, $j\in\Z_3$.
\subsection{Functional relations}
\label{sec:functional-relations}
We recall that the function $(W_j)_{j\in\Z_3}$  are defined in a domain $U$ satisfying the
assumptions from sections~\ref{ss:zeta} and~\ref{sec:asymptotic-solutions-1}.
\begin{Le} One has
  \begin{equation}
    \label{eq:three-W}
    W_0(z)+W_1(z)+W_2(z)=0, \qquad \forall z\in U.
  \end{equation}
\end{Le}
\begin{proof}
  Formula~(\ref{Airy}) and the definitions of the integration paths
  $(\gamma_j)_{j\in\Z_3}$ (see Fig.~\ref{fig:Airy}) imply that
  \begin{equation}
    \label{eq:three-Airy}
    w_0(\zeta)+w_1(\zeta)+w_2(\zeta)=0, \qquad \zeta\in\C.
  \end{equation}
  As the function $\zeta$ and all the coefficients
  $(A_l)_{l\in\N\cup\{0\}}$ and $(B_l)_{l\in\N\cup\{0\}}$ in
  representations~(\ref{w-sans-arg})--~(\ref{as-sol}) are independent
  of the choice of $w$, the solution of the Airy equation in this
  representation, the relation~\eqref{eq:three-Airy}
  implies~\eqref{eq:three-W}.
\end{proof}
\noindent Relation~(\ref{eq:three-W}) implies that
\begin{equation}
  \label{eq:three-Wron}
  (W_0(z),W_1(z))=(W_1(z),W_2(z))=(W_2(z),W_0(z)), \qquad \forall
  \{z,z+h\}\subset U,
\end{equation}
where $(f(z),\,g(z))=f(z+h)g(z)-g(z+h)f(z)$ is the difference Wronskian of $f$ and $g$.
\subsection{Estimates of $W_j$}
\label{sec:estimates-w_j}
To prove the existence of analytic solutions that admit asymptotic
expansions of the form~\eqref{eq:formal-sol}, we need rough estimates
of $(W_j)_{j\in\Z_3}$ in $U$. Therefore, we first introduce some
tools.
\subsubsection{Geometry}
\label{sec:geometry}
\begin{figure}
  \centering
  \includegraphics[height=5cm]{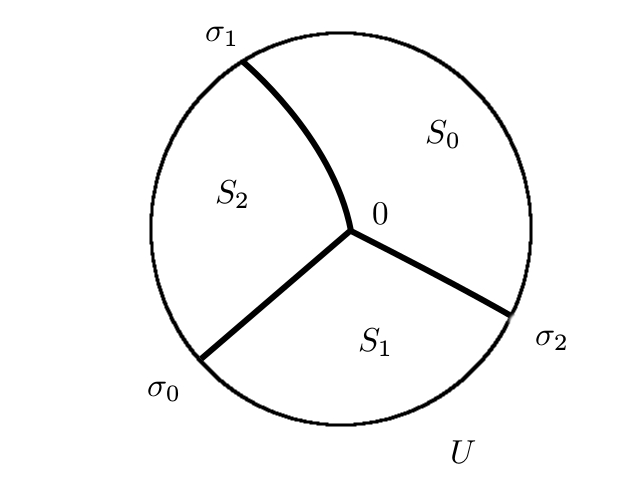}
  \caption{Stokes lines and sectors}\label{fig:Stokes}
\end{figure}

We recall that the function $\zeta$ defined in~\eqref{zeta} is
analytic in $U$ and bijectively maps $U$ onto $V=\zeta(U)$,
$\zeta(z_0)=0$.\\
We put
\begin{equation}
  \label{sigma:construction}
  \sigma_j=\zeta^{-1}(V\cap a_j), \quad a_j=e^{-2\pi i j/3}\,\R_-,\quad
  j\in\Z_3, 
\end{equation} 
where $\R_-=(-\infty,0]$.  The curves $(\sigma_j)_{j\in\Z_3}$ are
analytic. They all begin at $z_0$.  Any two of them do not intersect
except at $z_0$. The angles between these curves at $z_0$ are equal to
$2\pi/3$.\\
The curves $(\sigma_j)_{j\in\Z_3}$ cut the domain $U$ (a neighborhood
of $z_0$) into three simply connected subdomains that we call {\it
  sectors}. We denote them by $S_0$, $S_1$ and $S_2$ so that the
sector $S_0$ is bounded by $\sigma_1$ and $\sigma_2$, \ $S_1$ is
bounded by $\sigma_2$ and $\sigma_0$, and $S_2$ is bounded by
$\sigma_0$ and $\sigma_1$, see Fig.~\ref{fig:Stokes}.
%
%
Let
\begin{equation}
  \label{eq:U-j}
  U_j=U\setminus \sigma_j,\quad j\in\Z_3.
\end{equation}
These domains do not contain branch points of the complex momentum
$p$: the only branch point of $p$ in $U$ is $z=z_0$.  We shall use
\begin{Le}
  \label{le:actions-signs}
  For $j\in\Z_3$, there exists a branch $p_j$ of the complex momentum
  that is analytic in $U_j$ and such that $p_j(z_0)=0$ and
  \begin{enumerate}
  \item $\im \int_{z_0}^zp_j(z)\,dz >0$ inside $S_j$;
  \item $\im \int_{z_0}^zp_j(z)\,dz <0$ inside the two other
    sectors;
  \item $\im \int_{z_0}^zp_j(z)\,dz =0$ along the curves
    $\sigma_1$, $\sigma_2$ and $\sigma_0$ (in the case of $\sigma_j$,
    we mean the boundary values);
  \end{enumerate}
  Moreover, one has
  \begin{equation}
    \label{eq:p123}
    p_1=-p_0 \text{ in } \sigma_0\cup S_1\cup\sigma_2\cup
    S_0\cup\sigma_1, \quad  
    p_2=-p_0 \text{ in } \sigma_0\cup S_2\cup\sigma_1\cup
    S_0\cup\sigma_2. 
  \end{equation}
\end{Le}
In the WKB method, the curves $\sigma_j$, $j\in\Z_3$, are called {\it Stokes lines.}
\begin{proof} Let us check the first three points of
  Lemma~\ref{le:actions-signs} for $j=0$.  We recall that $\zeta$ is
  an analytic branch of the function~(\ref{zeta}).  We can assume that
  in~(\ref{zeta}) \ $p$ is a branch of the complex momentum analytic
  in $U_0$ and such that $p(z_0)=0$.
  \\
  Formulas~\eqref{sigma:construction} and the definition of $\zeta$
  imply that $\im \int_{z_0}^zp(z)\,dz = 0$ on any of the Stokes
  lines.  We note that
  \begin{equation}
    \label{eq:zeta-in-S-j}
    \zeta(S_j)=\{v\in V\;:\;v\ne 0,\; \arg v\in -2\pi  j/3+
    (-\pi/3,\pi/3)\},\quad j\in\Z_3. 
  \end{equation}
  This and the definition of $\zeta$ imply that
  $\im \int_{z_0}^zp(z)\,dz \ne 0$ in each of the sectors.\\
  In view of the analysis made in
  section~\ref{sec:analyt-branch-compl}, in $U_0$, we can choose an
  analytic branch $p_0$ of the complex momentum so that
  $\im \int_{z_0}^zp_0(z)\,dz>0$ in $S_0$. For $p_0$, the statements
  1. and 3. of Lemma~\ref{le:actions-signs} are obviously valid.\\
  To prove point 2., it suffices to check that
  $\im \int_{z_0}^zp_0(z)\,dz <0$ in the sectors $S_1$ and
  $S_2$. Therefore, we note that as $z\ne z_0$, $z\sim z_0$, crosses $\sigma_2$ moving
  from $S_0$ to $S_1$ the argument of $\zeta(z)$ decreases ($\zeta$
  vanishes only at $z=z_0$) as does the argument of
  $\int_{z_0}^zp_0(z)\,dz$. Therefore, point 2. of
  Lemma~\ref{le:actions-signs} follows from points 1. and 3.\\
  To complete the proof of Lemma~\ref{le:actions-signs}, we
  choose $p_1$ in the following way. First, we restrict $p_0$ to
  $S_1$; then, in $S_1$ we choose $p_1=-p_0$ and  continue $p_1$
  analytically from $S_1$ to $U_1$.  For the thus chosen $p_1$, we
  have
  \begin{equation*}
    \im \int_{z_0}^zp_1(z)\,dz=-\im \int_{z_0}^zp_0(z)\,dz
    >0,\quad z\in S_1.    
  \end{equation*}
  This proves point 1 for $p_1$. Point 2  and point 3 for
  $p_1$ are proved as for $p_0$ and $p$.\\
  To choose $p_2$, first, we restrict $p_0$ to $S_2$, then, in $S_2$
  we choose $p_2=-p_0$ and $p_2$ analytically from $S_2$ to $U_2$. To
  complete the proof of Lemma~\ref{le:actions-signs} for $p_2$, we
  reason as for $p_1$. We omit further details.
\end{proof}
\subsubsection{Estimates}
\label{as-sol:est}
For $j\in\Z_3$ and $z\in U_j$, we set
\begin{equation}
  \label{eq:Aj}
  \rho_j(z)=e^{\frac{i}h\int_{z_0}^zp_j(z')\,dz'}.
\end{equation}
We note that $\rho_j$ is continuous up to the cut along $\sigma_j$,
and the boundary values of its absolute value $|\rho_j|$ on both the
sides of the cut equal one. So, below, we consider $|\rho_j|$ as a
continuous function in $U$.

Let us recall that $H$ is defined by~(\ref{eq:H-def}). We set 
\begin{equation}\label{eq:delta}
\delta_j(z)=[H(W_j)](z),\quad z\in U,\qquad j\in Z_3.
\end{equation} 
\begin{Pro}
  \label{pro:W-est}
  For each $j\in\Z_3$, one has
  \begin{gather}
    \label{eq:W-est}
    |W_j(z)|\le C h^{1/3}|\rho_j(z)|,\quad  z\in U,\\
    \label{eq:delta-est}
    |\delta_j(z)|\le C
    h^{L+2+1/3}\,|\rho_j(z)|,\quad\{z,z+h,z-h\}\subset U,
  \end{gather}
  where $L$ is the order entering the definition of $W_j$,
  see~(\ref{as-sol}).
\end{Pro}
\noindent Proposition~\ref{pro:W-est} immediately follows from
formulas~(\ref{w-sans-arg})--(\ref{eq:as-sol}) with $w=w_j$ and
\begin{Le}
  \label{le:w:est} Let $j\in\Z_3$. Then one has
  \begin{equation}
    \label{eq:w-est}
    |w_j(h^{-\frac23}\zeta(z))|\le C|\rho_j(z)|,\quad   
    |w_j'(h^{-\frac23}\zeta(z))|\le Ch^{-\frac16}|\rho_j(z)|,\qquad  z\in U.
  \end{equation}
\end{Le}
\begin{proof}
  We prove~\eqref{eq:w-est} only for $j=0$. The other cases are
  treated similarly.
  We recall that $w_0={\rm Ai}$, that $\zeta$
  bijectively maps $U$ onto its image and that $\zeta(z_0)=0$
  (see~\eqref{zeta}). Clearly,
  \begin{equation}
    \label{eq:w-est:1}
    w_0(h^{-\frac23}\zeta(z))=O(1)\text{ and }
    w_0'(h^{-\frac23}\zeta(z))=O(1) \text{ if } |\zeta(z)|\le h^{\frac23}.
  \end{equation}
  Now we turn to the case where $|\zeta(z)|\ge h^{\frac23}$.  It suffices to 
  prove~\eqref{eq:w-est} in $U_0$.\\
  The asymptotic formulas~\eqref{ai:as:simple} imply that, 
  for $Z\in\{Z\in\C\setminus\R_-\,:\, |Z|\ge 1\}$,  one has 
  \begin{equation}
    \label{w-est}
    |w_0(Z)|\le C |Z|^{-\frac14}\;\left|e^{-\frac23\,Z^{\frac32}}\right|
    \text{ \ and \ } 
    |w_0'(Z)|\le C |Z|^{\frac14}\;\left|e^{-\frac23\, Z^{\frac32}}\right|,
  \end{equation} 
  where the $Z\mapsto Z^{\frac32}$  is analytic in $\C\setminus \R_-$ 
  and positive when $Z>0$.\\
  Estimate~\eqref{w-est} and the definition of $U_0$,
  see~(\ref{eq:U-j}), imply that, for $z\in U_0$ such that
  $|\zeta(z)|\ge h^{\frac23}$, one has
  \begin{equation}
    \label{eq:w-est:2}
    |w_0(h^{-\frac23}\zeta(z))|\le C\left|e^{-\frac2{3h}\,
        \zeta(z)^{\frac32}}\right| 
    \text{ and } 
    |w_0'(h^{-\frac23}\zeta(z))|\le
    Ch^{-\frac16}\left|e^{-\frac2{3h}\, \zeta(z)^{\frac32}}\right|,
  \end{equation}
  where $z\to \zeta(z)^{3/2}$ is analytic in $U_0$ and positive along
  $\alpha_0=\zeta^{-1}((0,\infty))$. 
  \\
  In view of the analysis made in
  section~\ref{sec:analyt-branch-compl},
  \begin{equation*}
    \zeta(z)^{\frac32}=\pm \frac{3i}2\int_{z_0}^zp_0(z')\,dz', \quad z\in  U_0.
  \end{equation*}
  As $\alpha_0=\zeta^{-1}((0,\infty))\subset S_0$, along $\alpha_0$ one has
  $\im \int_{z_0}^zp_0(z')\,dz'>0$.
  Therefore, in $U_0$ \
  $\zeta(z)^{\frac32}=-\frac{3i}2\int_{z_0}^zp_0(z')\,dz'$, and
  $\left|e^{-\frac2{3h}\,
      \zeta(z)^{\frac32}}\right|=|\rho_0(z)|$. This
  and~\eqref{eq:w-est:2} imply~\eqref{eq:w-est} for
  $ |\zeta(z)|\ge h^{-2/3}$. This and~\eqref{eq:w-est:1} imply the
  statement of the lemma.
\end{proof}
\subsection{Wronskians}
\label{sec:wronskians}
Below $\mathcal C\subset U$ is a closed disk independent of $h$ with the center at $z_0$. 
We now prove
\begin{Le}
  \label{le:2}
  For $\{z,z+h\}\subset \mathcal C$, as $h\to0$ one has
  \begin{equation*}
    (W_0(z),W_1(z))=h (w_0'(z)w_1(z)-w_0'(z)w_1(z))+O(h^{\frac53}).
  \end{equation*}
\end{Le}
\noindent Before proving Lemma~\ref{le:2}, we check
\begin{Le}
  \label{le:3}
  Let $j\in\Z_3$, and let  $w=w_j$. For $\{z,z+h\}\subset \mathcal C$, as $h\to0$ one has
  \begin{equation*}
    h^{\frac13}w_h\left|_{z+h}=h^{\frac13}\ch(\sqrt{\zeta(z)}\zeta'(z))\,w_h+
    g(z) h^{\frac23}\,w_h'+O(h^{\frac43}w_h)+O(h^{\frac53}w_h')\right.,
  \end{equation*}
  where we use the notations from~\eqref{w-sans-arg} and $g$ is
  defined in~\eqref{g:def}.
\end{Le}
\begin{proof}[Proof of Lemma~\ref{le:3}]
  We proceed as in the proof of
  Proposition~\ref{pro:Ahw}. Thus, we omit some details and
  concentrate on the new computations.\\
  Let $\gamma=\gamma_j$ be the integration path in~\eqref{Airy}. 
  Also, below we assume that  $z\in \mathcal C$ and that $h$ is sufficiently small.  
  The proof is broken into several steps.
  \\
  {\bf 1.} \ We compute
  \begin{equation}
    \label{Ew:int}
    \left.h^{\frac13}w_h \right|_{z+h}=
    \int_{\gamma}e^{-\frac1h\left(\frac{t^3}3-t\zeta(z)\right)}E(t,z,h)\,dt,\quad
    E(t,z,h)=e^{\frac{t}h(\zeta(z+h)-\zeta(z))}.
  \end{equation}
  \\
  {\bf 2.} \ We represent $E(t,z,h)$ in the form
  \begin{equation}\label{E-phi}
    E(t,z,h)=\alpha(z,h)+\beta(z,h)t+(t^2-\zeta(z)) \phi(t,z,h)
  \end{equation}
  with
  \begin{gather*}
    \alpha(z,h)=\frac{E(\sqrt{\zeta(z)},z,h)+E(-\sqrt{\zeta(z)},z,h)}2=
    \ch\left(\sqrt{\zeta}\frac{\zeta(z+h)-\zeta(z)}h\right),\\
    \beta(z,h)=\frac{E(\sqrt{\zeta(z)},z,h)-E(-\sqrt{\zeta(z)},z,h)}
    {2\sqrt{\zeta(z)}}= \frac1{\sqrt{\zeta}}\
    \sh\left(\sqrt{\zeta}\frac{\zeta(z+h)-\zeta(z)}h\right).
  \end{gather*}
  {\bf 3.} \ Clearly,
  \begin{equation}
    \label{as:alpha,beta}
    \alpha(z,h)=\ch(\sqrt{\zeta}\zeta')+O(h),\quad 
    \beta(z,h)=g(z)+O(h).
  \end{equation}
  \\
  {\bf 4.} \ Substituting~\eqref{E-phi} into~\eqref{Ew:int} and
  integrating by parts, we get
  \begin{equation}
    \label{int-of-phi1}
    \begin{split}
      \left.h^{\frac13}w_h \right|_{z+h}&=\alpha h^{\frac13}w_h +\beta
      h^{\frac23}w_h'
      +h\int_\gamma e^{-\frac1h\left(\frac{t^3}3-t\zeta(z)\right)} E_1(t,z,h)\,dt,\\
      &E_1(t,z,h)=\frac{\partial \phi}{\partial t}(t,z,h).
    \end{split}
  \end{equation}
  Reasoning as when proving Proposition~\ref{pro:Ahw}, we check that
  the last term in the right hand side of~\eqref{int-of-phi1} is
  $O(h^{\frac43}w)+O(h^{\frac53}w')$. Lemma~\ref{le:3} follows from
  this estimate, asymptotics~\eqref{as:alpha,beta} and
  representation~\eqref{int-of-phi1}.
\end{proof}
\noindent Now, we turn to the proof of Lemma~\ref{le:2}.
\begin{proof}[Proof of Lemma~\ref{le:2}]
  Below we assume that $\{z,z+h\}\subset \mathcal C$ and that $h$ is 
  sufficiently small. Using~(\ref{as-sol}) and Lemma~\ref{le:3}, we compute
  \begin{equation*}
    \begin{split}
      (W_0(z),&W_1(z))=A_0^2(z)\\
      &\times\left(\,( h^{\frac13}\ch(\sqrt{\zeta}\zeta')\,w_0
        +gh^{\frac23}\,w_0'+O(h^{\frac43}w_0)+O(h^{\frac53}w_0'))\right.\\
      &\hspace{4cm}\cdot(h^{\frac13}w_1+O(h^{\frac43}w_1)
      +O(h^{\frac53}w_1'))\\
      &\hspace{1.5cm}-(h^{\frac13}\ch(\sqrt{\zeta}\zeta')\,w_1
      +g h^{\frac23}\,w_1'+O(h^{\frac43}w_1)+O(h^{\frac53}w_1'))\\
      &\hspace{5.5cm}\left.\cdot(h^{\frac13}w_0+O(h^{\frac43}w_0)
        +O(h^{\frac53}w_0'))\right).
    \end{split}
  \end{equation*}
  Here, $w_j=w_j(h^{-\frac23}\zeta(z))$, $j\in\Z_3$, $\zeta=\zeta(z)$,
  and $g=g(z)$.\\
  Now, we assume that
  $z\in \sigma_0\cup S_1\cup\sigma_2\cup S_0\cup \sigma_1$.  Then,
  by~(\ref{eq:w-est}) and~(\ref{eq:p123})
  \begin{equation}
    |w_0w_1|\le C,\quad |w_0'w_1|\le C  h^{-\frac16},\quad |w_0w_1'|\le C h^{-\frac16}, 
   \quad |w_0'w_1'|\le C h^{-\frac13},
  \end{equation}
  and we get
  \begin{equation*}
    (W_0(z),W_1(z))=hg(z)A_0^2(z) ( w_0'w_1-w_0w_1') +O(h^{\frac53}).
  \end{equation*}
  In view of~(\ref{eq:A0}), Lemma~\ref{le:2} is proved for
  $z\in \sigma_0\cup S_1\cup\sigma_2\cup S_0\cup \sigma_1$.\\
  When $z\in \sigma_0\cup S_2\cup\sigma_1\cup S_0\cup \sigma_2$, we
  similarly get
  \begin{equation}\label{wrosk:W0W2}
    (W_0(z),W_2(z))=hg(z)A_0^2(z) ( w_0'w_2-w'_2w_0) +O(h^{\frac53}).
  \end{equation}
  In view of~\eqref{eq:three-Wron}, $(W_0,W_2)=-(W_0,W_1)$ and
  relation~(\ref{eq:three-Airy}) imply that
  $w_0'w_2-w_0w_2'=-(w_0'w_1-w_0w_1')$. Therefore, Lemma~\ref{le:2}
  for $z\in S_2$ follows from~\eqref{wrosk:W0W2}.  This completes the
  proof of Lemma~\ref{le:2}.
\end{proof}
\section{Solutions to~\eqref{eq:main} on precanonical domains}
\label{s:CD}
Fix $L\in \N$. Here we construct solutions $(\psi_j)_{j\in\Z_3}$ to
equation~\eqref{eq:main} that, up to $O(h^L)$, coincide with
$(W_j)_{j\in\Z_3}$, the functions from Theorem~\ref{th:as-sol}.  The
result of this section is preliminary: we only construct the
$(\psi_j)_{j\in\Z_3}$ on some subdomains of $U$.
\subsection{The result of this section}
\label{sec:result-this-section}
\subsubsection{Notations and some definitions}
\label{sec:notat-some-defin}
First, to formulate the results of this section, we introduce some
notations and recall some definitions related to the complex WKB
method for difference equations, see, for example,~\cite{F-Shch:15}.\\
For $z\in \C$, we let $x=\re z$ and $y=\im z$.\\
A curve $\gamma\subset \C$ is called {\it vertical}, if $z$ is
a piecewise continuously differentiable function of $y$ along
$\gamma$.\\
Let $\gamma\subset U$ be a regular vertical curve parameterized by $z=z(y)$. Let
$p$ be a branch of the complex momentum continuous on $\gamma$. 
The curve $\gamma$ is {\it precanonical} with respect to $p$, 
if the function $ y\mapsto \im\int_{z_0}^{z(y)}p(z)\,dz$ is non decreasing and the
function $ y\mapsto \im\int_{z_0}^{z(y)}(p(z)-\pi)\,dz$ is
non increasing.\\
Let $d>0$. For $M\subset C$ we define the {\it horizontal
  $d$-neighborhood} of $M$ to be the set $M^d:=M+[-d,d]$ and
$M^{-d}:=(M^d-d)\cap M\cap(M^d+d)$.\\
We recall that, for $j\in\Z_3$, the sector $S_j$ and the Stokes line
$\sigma_j$ are shown in Fig.~\ref{fig:Stokes}. For $j\in\Z_3$, we
denote by $S_{j,j+1}$ the closure of the domain $S_j\cup S_{j+1}$
without the boundary of $U$. For example, one has
\begin{equation*}
  S_{1,2}=\sigma_1\cup S_2\cup \sigma_0\cup S_1 \cup \sigma_2.
\end{equation*}
We also note that relations~\eqref{eq:p123} imply that
\begin{equation}
  \label{eq:rho123}
  |\rho_j(z)\rho_{j+1}(z)|=1 ,\qquad z\in S_{j,j+1},\quad j\in\Z_3.
\end{equation}
Let $r_1<r_2$.  We set $S(r_1,r_2)= \{z\in \C\,:\, r_1<\im z<r_2\}$.
\subsubsection{The main result of the section}
\label{sec:main-result-section}
One has
\begin{Th}
  \label{th:CD}
  Let $j\in\Z_3$, $L\in\N$, $c\in (1,2)$ and $r>0$.
  Let $K\subset S_{j,j+1}$ be a regular simply connected domain bounded by
  two curves having common endpoints $z_1$ and $z_2$ and both
  precanonical with
  respect to either the branch $p_j$ or $p_{j+1}$.\\
  Then, for sufficiently small $h$, there exist two solutions $\psi_j$ and
  $\psi_{j+1}$ to~\eqref{eq:main} that are analytic in $K^{ch}$ and
  that, in $K^{ch}\cap S(\im z_1+rh,\,\im z_2-rh)$ admit the
  asymptotic representations
  \begin{equation}
    \label{eq:sol-as-exp}
    \psi_l(z)=W_l(z)+O(|\rho_l|h^{L+1+\frac13}),\quad l\in\{j,j+1\},
  \end{equation}
  where $W_l$ is the function described in Theorem~\ref{th:as-sol}
  and corresponding to $w_l$ and  the order $L$.
\end{Th}
\noindent Let us discuss the solutions $\psi_j$ and $\psi_{j+1}$
described in Theorem~\ref{th:CD}.
\begin{Cor}
  \label{cor:th:CD}
  In the case of Theorem~\ref{th:CD}, the solutions $\psi_j$ and
  $\psi_{j+1}$ can be analytically continued to
  $U\cap S(\im z_1,\,\im z_2)$. Let $r>0$. As $h\to 0$, one has
  \begin{equation}\label{psi-Wrons}
    (\psi_j(z),\psi_{j+1}(z))=(W_j(z),W_{j+1}(z))+O(h^{L+1+\frac23})
  \end{equation}
  in  $K^{ch}\cap S(\im z_1+rh,\,\im z_2-rh)$.
%
%
%
\end{Cor}

\begin{proof}
  The solutions being analytic in $K^{ch}$ with $c>1$, they can be
  analytically continued to $U\cap S(\im z_1,\,\im z_2)$ just be
  means of   equation~\eqref{eq:main}.\\
%
%
  We fix $l\in\Z_3$ and note that, for all $z$ in a compact set
  $\mathcal C\subset U$, for sufficiently small $h$, one has
  $|\rho_l(z+h)|/|\rho_l(z)|\le C$.  For
  $z\in K^{ch}\cap S(\im z_1+rh,\,\im z_2-rh)$,
  representation~\eqref{psi-Wrons} follows from this observation and
  from~(\ref{eq:W-est}),~(\ref{eq:sol-as-exp})
  and~(\ref{eq:rho123}).
\end{proof}
The remainder of this section is devoted to the proof of
Theorem~\ref{th:CD}.\\
For the sake of definiteness, when proving Theorem~\ref{th:CD}, we
assume that $j=0$ and that the two curves from Theorem~\ref{th:CD} are
precanonical with to respect to the branch $p_0$. The other cases are
treated similarly.\\
Below, $K$ is as in the theorem (for $j=0$),  it is bounded by the 
precanonical curves $\gamma_1$ and
$\gamma_2$, and their common endpoints satisfy $\im z_1<\im z_2$. 
Finally, $h$ is supposed to be sufficiently small.
\subsection{Ideas of the proof}
\label{sec:ideas-proof}
In the present section, we describe the construction of the solution
$\psi_0$. The solution  $\psi_1$ is constructed similarly.\\
Let us assume that $\psi_0$ is a solution to~\eqref{eq:main} analytic
in $K^{ch}$ that we expect to be close to $W_0$.  Let us recall that
$\delta_0=HW_0$.  Clearly, $ \Delta_0:=W_0-\psi_0$ satisfies the
equation
\begin{equation}
  \label{eq:HDelta}
  H(\Delta_0)(z)=\delta_0(z), \quad \{z-h,\,z,\,z+h\}\subset K^{ch}.  
\end{equation}
For $z\in K^{ch}$, let $\gamma(z)$ denote a vertical curve in $K^{ch}$
that contains $z$ and connects $z_1$ and $z_2$.
We construct a solution to the equation for $\Delta_0$ in the form
\begin{equation}
  \label{eq:ansatz}
  \Delta_0=R_0\,g_0 \quad \text{where \ } 
  R_0g_0\,(z):=\int_{\gamma(z)}r_0(z,\zeta)g_0(\zeta)\,d\zeta,
\end{equation}
\begin{equation}
  \label{eq:r}
  r_0(z,\zeta)=\frac1{2ih}\,\frac{W_0(z)W_1(\zeta)-W_0(\zeta)W_1(z)}
  {(W_0(\zeta),\,W_1(\zeta))}\,
  \theta_0\left(\frac{\zeta-z}h\right),\quad 
  \theta_0(t)=\cot(\pi t)- i.
\end{equation}
Here, $(W_0(\zeta),\,W_1(\zeta))$ is the difference Wronskian of $W_0$ and $W_1$.
The choice of Ansatz~\eqref{eq:ansatz} is explained by
\begin{Le}
  \label{le:parametrix}
  Let $0<\beta<1$.  Let $f$ be a function defined and analytic in
  $U\cap S(\im z_1,\im z_2)$ and such that the expression
  \begin{equation}\label{eq:fa}
    f_\beta(z)=(z-z_1)^\beta (z-z_2)^\beta f(z)
  \end{equation}
  is bounded. Then, if $\{z-h,z,z+h\} \subset U\cap S(\im z_1,\im z_2)$, one has
  \begin{equation}
    \label{eq:parametrix}
    HR_0f(z)=f(z)+D_0f(z), \quad
    D_0f(z)=\int_{\gamma(z)}d_0(z,\zeta)f(\zeta)\,d\zeta
  \end{equation}
  where
  \begin{equation}
    \label{eq:d}
    d_0(z,\zeta)=\frac1{2ih}
    \frac{\delta_0(z)W_1(\zeta)-W_0(\zeta)\delta_1(z)}
    {(W_0(\zeta),\,W_1(\zeta))}\,
    \theta_0\left(\frac{\zeta-z-0}h\right),
  \end{equation}
  and $\delta_j:=HW_j$ are the ``error'' terms estimated
  in~(\ref{eq:delta-est}).  The function $D_0f$ is analytic in
  $U\cap S(\im z_1,\im z_2)$.
\end{Le}
\begin{proof} The analyticity of $f$ and the boundedness of $f_\beta$
  imply that $D_0f$ is well defined and analytic in
  $U\cap S(\im z_1,\im z_2)$.  The relation $HR_0f=f+D_0f$ follows
  from the residue theorem. We omit further details.
\end{proof}
 We note that an operator similar to $R_0$ was introduced
in~\cite{F-Sh:2017}, but, was not studied for small $h$.\\
In view of Lemma~\ref{le:parametrix} and the formulas
$\Delta_0=R_0g_0$ and $H(\Delta_0)=\delta_0$, we can expect that in
$K^{ch}$
\begin{equation}
  \label{eq:eq:g0}
  g_0+D_0g_0=\delta_0.
\end{equation}
Roughly, to prove Theorem~\ref{th:CD}, we consider~\eqref{eq:eq:g0} as
an equation on a vertical curve $\gamma$. It appears that if $\gamma$
is precanonical, the operator $D_0$ is small.  This enables us to
construct a solution $\psi_0$ to equation~\eqref{eq:eq:g0} on
$\gamma_1$.  Next, we check that it is analytic in $K^{ch}$,
satisfies~\eqref{eq:main} and admits the asymptotic
representation~\eqref{eq:sol-as-exp}. The solution $\psi_1$ is
constructed similarly.
\subsection{The integral operator $D_0$}
\label{sec:integral-operator}
The aim of this section is to estimate the ope\-rator norm of $D_0$ in a
suitable functional space.\\
Let $\gamma$ be either $\gamma_1$ or $\gamma_2$.  We fix $\alpha\in(0,1)$ and
define the strip
\begin {equation*}
  \Pi_{\gamma, \alpha} = \gamma\setminus\{z_1,z_2\}+[-\alpha h,\alpha h].
\end {equation*}
We recall that $z\to |\rho_0(z)|$ defined in $U_0$ is a continuous
function in $U$.  We fix $ 0 <\beta <1 $ and let
$H_{\gamma, \alpha, \beta}$ be the linear space of functions analytic in
$\Pi_{\gamma, \alpha}$ and having finite norm
\begin{equation}
  \label{def:norm}
  \| f \| = \sup_ {z \in \Pi_ {\gamma,  \alpha}}\frac{|f_\beta(z) |}{|\rho_0(z)|}
\end{equation}
$f_\beta$ being defined in~(\ref{eq:fa}).
Obviously, endowed with this norm, $H_{\gamma, \alpha, \beta}$ is a Banach
space.\\
For $f\in H_{\gamma, \alpha,\beta}$, we define $D_0f$ by the formula
in~\eqref{eq:parametrix}, where $\gamma(z)$ is a vertical curve that
connects the points $z_1$ and $z_2$ in $\Pi_{\gamma, \alpha}$ and passes
through $z$. The function $D_0f$ is then analytic in
$\Pi_{\gamma, \alpha}$. One has
\begin {Pro}
  \label{pro:AB}
  For sufficiently small $h$
  \begin{equation*}
    \| D_0 \|_{H _ {\gamma,  \alpha, \beta} \to H _ {\gamma,  \alpha, \beta}} \le Ch^{L+\frac23}.
  \end{equation*}
\end {Pro}
\noindent
The remainder of this subsection is devoted to the proof of
Proposition~\ref{pro:AB}. Therefore, for $f\in H_{\gamma,\alpha,\beta}$, we
estimate $Rf(z)$.  Up to the end of this subsection, we assume that
$\{z, \zeta\}\subset\Pi_{\gamma,\alpha}$ and that $h$ is sufficiently
small.
\subsubsection{Auxiliary lemma}
\label{sec:auxiliary-lemma}
When estimating $Rf(z)$, we use
\begin{Le}
  \label{le:1}
  For $q>0$, there exists $C>0$ such that
  \begin{equation}
    \label{est:d0}
    \sup_{\substack{\D\{z, \zeta\}\subset\Pi_{\gamma,\alpha} \\
        \D\min_{k\in\Z}|\zeta-z-kh|\geq
        qh}}\left|\frac{\rho_0(\zeta)}{\rho_0(z)}\,
      d_0(z,\zeta)\right|\leq C h^{L+\frac23}.
  \end{equation}
\end{Le}
\begin{proof}
  We proceed in several steps.\\
  {\bf 1.} \ Proposition~\ref{pro:W-est} and Lemma~\ref{le:2} imply
  that
  \begin{equation}\label{aux:1}
    \left|\frac{\rho_0(\zeta)}{\rho_0(z)}\, d_0(z,\zeta)\right|\le C
    h^{L+\frac23}
    \left(|\rho_1(\zeta)\rho_0(\zeta)|+
      \frac{|\rho_1(z)\rho_0^2(\zeta)|}{|\rho_0(z)|}\right) 
    \left|\theta_0\left(\frac{\zeta-z-0}h\right)\right|.
  \end{equation}
  {\bf 2.} \ Recall that in $ S_{0,1}$ one has $\rho_0(z)\rho_1(z)=1$ 
  (see~(\ref{eq:rho123})).  As $\gamma_1,\,\gamma_2\subset S_{0,1}$,
  one has $|\rho_0(z)||\rho_1(z)|\le C$ for $z\in \Pi_{\gamma, \alpha}$.  Therefore,
  \begin{equation}
    \label{est:r:aux}
    \left|\frac{\rho_0(\zeta)}{\rho_0(z)}\, d_0(z,\zeta)\right|\le C
    h^{L+\frac23} \left(1+
      e(z,\zeta)\right)
    \left|\theta_0\left(\frac{\zeta-z-0}h\right)\right|,\quad    
    e(z,\zeta)= \left|\frac{\rho_0(\zeta)}{\rho_0(z)}\right|^2.
  \end{equation}
  {\bf 3.} \ For $z\in \Pi_{\gamma, \alpha}$, we define $z_\perp\in \gamma$ so
  that $\im z_\perp=\im z$.  We have
  \begin{equation*}
    |e(z,\zeta)|\le C\,
    \left|\exp\left(\frac{2i}h\int_{z_\perp}^{\zeta_\perp}
        p_0(z')\,dz'\right)\right|.
  \end{equation*}
  {\bf 4.} \ On the complex plane outside a fixed neighborhood of the
  points $\Z$, we have the estimate
  \begin {equation*}
    |\theta_0(z)|=| \cot (\pi z) -i | \le C \, \begin {cases} 1, & \im z \ge 0; \\
      e ^ {2 \pi \im z}, & \im z \le 0. \end {cases}
  \end {equation*}
  Therefore, for $\zeta$ outside the $(qh)$-neighborhood of $z+h\Z$,
  we get
  \begin {gather*}
    \left|\theta_0\left(\frac{\zeta-z-0}h\right)\right| \le C\quad
    \text{and}\\
    \left|e(z,\zeta)\, \theta_0\left(\frac{\zeta-z-0}h\right)\right|
    \le C \,
    \begin {cases} e ^ {- \frac {2} {h} \im \int \limits_{z_\perp}
        ^{\zeta_\perp} p \, dz} &\text{ if }\im (\zeta- z)\ge 0; \\
      e ^ {\frac {2} {h} \im \int \limits_{\zeta_\perp}^{z_\perp} (p-
        \pi) \, dz} &\text{ if } \im(\zeta-z)\le 0. \end {cases}
  \end {gather*}
  {\bf 5.} \ As $\gamma$ is a precanonical curve, we finally get
  \begin {equation}\label{est:theta}
    \left|\theta_0\left(\frac{\zeta-z-0}h\right)\right| \le
    C,\quad\text{and}\quad 
    \left|e(z,\zeta)\,
      \theta_0\left(\frac{\zeta-z-0}h\right)\right| \le C.
  \end {equation}
  This and~\eqref{est:r:aux} imply~\eqref{est:d0}.
\end{proof}
\subsubsection{Estimates in the strip
  $S(\im z_1+h/2,\, \im z_2-h/2)$}
\label{sss:est:D0:1}
When $z\in S(\im z_1+h/2,\,\im z_2-h/2)$, we prove
\begin{equation}
  \label{est:Af} 
  |\rho_0(z)^{-1}D_0f(z)|\leq Ch^{L+\frac23}\|f\|.
\end{equation}
First, we assume that $z$ is between the curves $\gamma+\alpha h/2$ and
$\gamma+\alpha h$. Then, one can deform the integration path $\gamma(z)$
in~(\ref{eq:parametrix}) to $\gamma$.  The distance between the poles
of $d_0$ and $\gamma$ is larger than $Ch$. This,~\eqref{def:norm}
and~\eqref{est:d0} imply~\eqref{est:Af}.\\
Next, we assume that $z$ is either between the curves $\gamma$ and
$\gamma+ah/2$ or on one of them. In this case,
in~(\ref{eq:parametrix}) we can replace the integration path $\gamma$
by $\tilde \gamma$ where
\begin{itemize}
\item $\tilde \gamma$ is a continuous curve that connects $z_1$ to
  $z_2$,
\item $\tilde \gamma$ coincides with $\gamma-\alpha h/2$ in the strip
  $\{\im z_1+h/2\le \im z\le \im z_2-h/2\}$,
\item outside this strip, $\tilde \gamma$ consists of two segments of
  straight lines.
\end{itemize}
Reasoning as above on this new integral, we again
obtain~\eqref{est:Af}.\\
Let us assume now that $z$ is to the left of $\gamma$. We note that,
by the Residue theorem, the integral in~\eqref{eq:parametrix}
decomposes as the sum of
\begin{equation*}
  -\frac{\delta_0(z)W_1(z)-W_0(z)\delta_1(z)}{(W_0(z),\,W_1(z))}f(z)=
  O(h^{L+1+\frac23})\,f(z)
\end{equation*}
and the integral defined by~\eqref{eq:parametrix}--(\ref{eq:d}) with
$\theta( (\zeta-(z+0))/h)$ replaced with $\theta((\zeta-(z-0))/h)$.\\
This new integral for $z$ to the left of $\gamma$ is analyzed as
above. This completes the proof of~\eqref{est:Af}.
\subsubsection{Estimates in 
  $S(\im z_1,\,\im z_1+h/2)$ and
  $S(\im z_2-h/2,\, \im z_2)$}
\label{sec:estimates-domains-im}
Both domains are treated similarly. So, we detail only the analysis
for the first one. We prove that
\begin{equation}
  \label{eq:Af:end}
  |\rho_0(z)^{-1}(D_0f)_\beta(z)|\le  C \,h^{L+\frac23}\, \|f\|. 
\end{equation}
For $z$ between $\gamma+\alpha h/2$ and $\gamma+\alpha h$, reasoning as in
section~\ref{sss:est:D0:1}, one obtains~(\ref{est:Af}) that
implies~\eqref{eq:Af:end}.\\
For $z$ between $\gamma$ and $\gamma+\alpha h/2$, by contour deformation,
the integration path in~\eqref{eq:parametrix} is replaced with
$\tilde \gamma$ defined in section~\ref{sss:est:D0:1}. We, thus, write
$D_0f$ as the sum of an integral, say $A$, over the part of
$\tilde\gamma\cap\{\im \zeta\le z_1+h/2\}$ and an integral, say $B$,
over the part of $\tilde\gamma\cap\{\im z_1+h/2\le \im\zeta\}$.\\
Reasoning as in section~\ref{sss:est:D0:1}, we estimate $B$ and obtain
\begin{equation}
  \label{eq:1}
  |\rho_0(z)^{-1}B|\le Ch^{L+\frac23}\|f\|. 
\end{equation}
Let us turn to $A$. We again use~\eqref{aux:1}. 
Now, both $|z-z_1|$ and $|\zeta-z_1|$ are bounded
by $Ch$; thus, $ |\rho_0(z)/ \rho_0(\zeta)|\le C$.  Furthermore, for
such $z$, only one pole of the integrand, the pole at the point $z$,
can approach the integration path in $A$; the other poles stay at a
distance greater than $Ch$ from it. Therefore, we get
\begin {equation*}
  | \rho_0(z)^{-1}A| \le C \,h^{L+1+\frac23}\, \|f\| 
  \int_{\im z\le \im z_1+h/2} \frac {|d \zeta|}{|z- \zeta|\,
    |\zeta-z_1|^\beta},
\end {equation*}
where we integrate along $\tilde \gamma$. Changing variable
$t=(\zeta-z_1)/|z-z_1|$, one checks that the last integral is bounded
by $C/|z-z_1|^\beta$. Thus,
\begin{equation*}
  | \rho_0(z)^{-1}A | \le C \,h^{L+1+\frac23}\, \|f\| /|z-z_1|^\beta. 
\end{equation*}
This and~\eqref{eq:1} yields~\eqref{eq:Af:end}.

We omit further details and note only that, to prove~\eqref{eq:Af:end}
when $z$ is to the left of $\gamma$, we first transform the integral
from~(\ref{eq:parametrix}) as when doing the estimations in the strip
$S(\im z_1+h/2,\,\im z_2-h/2)$ (see the end of the
section~\ref{sss:est:D0:1}).
\subsubsection{Completing the proof of Proposition~\ref{pro:AB}}
\label{sec:compl-proof-prop}
Proposition~\ref{pro:AB} follows from estimates~\eqref{est:Af}
and~\eqref{eq:Af:end}.
\subsection{Solutions to the integral equation~(\ref{eq:eq:g0})}
\label{ss:g0}
Consider the integral equation~(\ref{eq:eq:g0}) in
$ H _ {\gamma, \alpha, \beta} $.  Proposition~\ref{pro:AB} and the estimate
for $ \delta_0$ from~(\ref{eq:delta-est}) imply
\begin {Le}
  \label{le:solution-of-integral-eq}
  For sufficiently small $ h $, the equation~\eqref{eq:eq:g0} has a
  unique solution $ g_0$ in $ H _ {\gamma, \alpha, \beta} $. It satisfies
  \begin {equation}
    \label{est:g_0}
    \|g_0(z)\| = O(h^{L+2+\frac13}).
  \end {equation}
\end {Le}
\noindent Moreover, one has
\begin{Le}
  \label{le:4}
  The solution $g_0$, constructed in
  Lemma~\ref{le:solution-of-integral-eq} for the curve
  $\gamma=\gamma_1$, can be analytically continued to the domain
  $K^{\alpha h}$. It then satisfies~(\ref{eq:eq:g0}) and in $K^{\alpha h}$
  \begin{equation}
    \label{est:g0}
    |(z_1-z)(z_2-z)|^\beta
    \frac{|g_0(z)|}{|\rho_0(z)|}\leq Ch^{L+2+\frac13}.
  \end{equation}
\end{Le}
\begin{proof}
  The proof is divided into four parts.\\
  {\bf 1.}\ As $g_0$ is analytic in $\Pi_{\gamma_1,\alpha}$, it suffices to
  continue it to the right of $\gamma_1$. The function
  $\zeta\to\theta_0\left(\frac{\zeta-z-0}h\right)$ has all its poles
  in $z+0+h\Z$. Hence, for $z$ between $\gamma_1$ and $\gamma_1+h$, we
  can define $D_0g_0$ by means of~\eqref{eq:parametrix} with
  $\gamma(z)=\gamma_1$, and $D_0g_0$ appears to be analytic between
  $\gamma_1$ and $\gamma_1+h$.\\
  As $g_0$ is analytic between $\gamma$ and $\gamma+\alpha h$, to define
  $D_0g_0$ for $z$ between $\gamma_1+\alpha h$ and $\gamma_1+(1+\alpha)h$, we
  can deform the path of the integral in~\eqref{eq:parametrix} to a
  vertical curve connecting $z_1$ to $z_2$ staying between $ \gamma_1$
  and $\gamma_1+\alpha h$.  Thus,~(\ref{eq:parametrix}) implies that
  $D_0g_0$ is analytic in $z$ between $\gamma_1$ and $\gamma_1+\alpha h+h$.
  In view of equation~(\ref{eq:eq:g0}), this implies that $g_0$ itself
  is analytic to the left of $\gamma+(\alpha+1)h$.  Reasoning in this way
  inductively, one shows that $g_0$ and $D_0g_0$ are analytic between
  $\gamma$ and $\gamma+(\alpha+2)h$, between $\gamma$ and $\gamma+(\alpha+3)h$
  and so on. As a result, one sees that $g_0$ and $D_0g_0$ are
  analytic in $K^{\alpha h}$ to the right of $\gamma$ and
  satisfy~(\ref{eq:eq:g0}) for all $z\in K^{\alpha h}$.
  \\
  {\bf 2.} \ Clearly, $g_0$ is analytic in $\Pi_{\gamma_2,\alpha}$, the
  expression $|(z_1-z)(z_2-z)|^\beta|g_0(z)|$ stays bounded in
  $\Pi_{\gamma_2,\alpha}$ (as $\gamma_1$ and $\gamma_2$ have common ends),
  and $g_0$ satisfies equation~\eqref{eq:eq:g0} along
  $\gamma=\gamma_2$. By Lemma~\ref{le:solution-of-integral-eq}, for
  sufficiently small $h$, this equation has a unique solution in
  $H_{\gamma_2,\alpha,\beta}$ which, thus, coincides with $g_0$. Hence,
  $g_0$ satisfies~\eqref{est:g_0} with the norm of
  $H_{\gamma_2,\alpha,\beta}$.
  \\
  {\bf 3.}  In view of the previous step, $g_0$ satisfies~(\ref{est:g0}) 
  in $\Pi_{\gamma_1,\alpha}\cup\Pi_{\gamma_2,\alpha}$. This and the maximum 
  principle for analytic functions imply that $g_0$ satisfies~(\ref{est:g0}) also in
  the domain bounded by $\gamma_1$ and $\gamma_2$, i.e., in $K$.\\
  The proof of Lemma~\ref{le:4} is complete.
\end{proof}
\subsection{The solution to the difference equation}
\label{sec:solut-diff-equat}
We define $\Delta_0$ by~\eqref{eq:ansatz} in terms of $g_0$
constructed in section~\ref{ss:g0}. One has
\begin{Le}
  \label{le:5}
  The function $\Delta_0$ can be analytically continued to
  $K^{(1+\alpha)h}$ where it satisfies equation~\eqref{eq:HDelta}.\\
  Let $0<c<1+\alpha$ and $r>0$. In $K^{ch}\cap S(\im z_1+rh,\im z_2-rh)$, one has
  \begin{equation}
    \label{est:Delta0}
    |\Delta_0(z)|\le C |\rho_0(z) h^{L+1}|.
  \end{equation}
\end{Le}
\begin{proof}
  By~(\ref{est:g0}) the function
  $z \mapsto|(z_1-z)(z_2-z)|^\beta|g_0(z)|$ is bounded in $K^{\alpha
    h}$. For a given $z$, the poles of the kernel
  in~\eqref{eq:ansatz} are contained in $z+h
  (\Z\setminus\{0\})$. Thus, the function
  $\Delta_0$ is analytic in $\left(K^{\alpha h}\right)^h=K^{(1+\alpha)h}$.\\
  By means of the Residue theorem, one checks
  that $H\Delta_0= HR_0g_0$ is equal to $g_0+D_0g_0$ if
  $z,z\pm h\in K^{(1+\alpha)h}$. As $g_0$ satisfies~(\ref{eq:eq:g0}) in
  $K^{\alpha h}$, we obtain~(\ref{eq:HDelta}) if
  $z,z\pm h\in K^{(1+\alpha)h}$.\\
  To prove~\eqref{est:Delta0}, we estimate $R_0g_0$ in the same way as
  in section~\ref{sss:est:D0:1} we estimated $D_0f$. So, we omit
  further details and only note that
  \begin{enumerate}
  \item outside $(Ch)$-neighborhood of the set $z+h\Z$, instead
    of~(\ref{est:d0}) we obtain
    \begin{equation*}
      \left|\frac{\rho_0(\zeta)}{\rho_0(z)}\, r_0(z,\zeta)\right|\le C
      h^{- \frac43};  
    \end{equation*}
  \item on the diagonal $\{\zeta=z\}$, $r_0$, the kernel of $R_0$, is
    analytic whereas $d_0$, the kernel of $D_0$, has a pole.  This
    simplifies the estimates of $(R_0g_0)(z)$ to the left of
    $\gamma_1$.
  \end{enumerate}
\end{proof}
\noindent Having constructed $\Delta_0$, we construct a solution
$\psi_0$ to equation~(\ref{eq:main}) setting $\psi_0=W_0-\Delta_0$,
see ~(\ref{eq:HDelta}). Let $c\in (0,2)$. In view of~(\ref{est:Delta0}), 
one has 
\begin{equation}\label{almost}
  \psi_0(z)=W_0(z)+O(|\rho_0(z)| h^{L+1}),\qquad z\in K^{ch}\cap S(\im
  z_1+rh,\,\im z_2-rh).
\end{equation}
In view of~(\ref{eq:w-est}), estimate~\eqref{almost}
implies~(\ref{eq:sol-as-exp}) with $L$ replaced with $L-1$.  As we
could choose a larger $L$, this actually completes the proof of the
statement of Theorem~\ref{th:CD} on the solution $\psi_0$.
\subsection{The second solution}
\label{sec:second-solution}
Mutatis mutandis, the construction of the solution $\psi_1$ repeats
that of $\psi_0$. We omit further details and mention only that, in
this case,
\begin{itemize}
\item we set $\psi_1=W_1-R_1g_1$ where $R_1$ is an integral operator
  with the kernel
  \begin{equation*}
    r_1(z,\zeta)=\frac1{2ih}\,\frac{W_0(z)W_1(\zeta)-W_0(\zeta)W_1(z)}
    {(W_0(\zeta),\,W_1(\zeta))}\,
    \theta_1\left(\frac{\zeta-z}h\right), \quad  \theta_1(t)=\cot(\pi t)+ i;
  \end{equation*}
\item instead of~\eqref{def:norm}, we use the norm
  $\D \| f\|=\sup_ {z \in \Pi_ {\gamma, \alpha}} \frac{| f_\beta(z) |}{|\rho_1(z)|}$.
\end{itemize}
\section{Proof of the main Theorem}
\label{s:proof:main}
In this section we finally prove Theorem~\ref{main:th}. 
We recall that in $U$ there are three Stokes lines beginning at
$z_0$. They are analytic curves, and the angle between any two of them
at $z_0$ is equal to $2\pi/3$. So, possibly  reducing $U$ somewhat, 
we can assume that at least two of them form a
vertical curve. We prove the theorem in the case where these are $\sigma_1$
and $\sigma_2$, and  $\sigma_1$ goes upwards from
$z_0$, i.e., the vector tangent to $\sigma_1$ at $z_0$ is
directed  in the upper half-plane. Mutatis mutandis, the other cases are treated 
in the same way. Moreover, we assume that the tangent vector to $\sigma_0$ is either 
directed in the lower half-plane or is parallel to the real line and directed to the left.
Then the curves $\sigma_j$, $j\in \Z_3$, correspond to Fig.~\ref{fig:Stokes}. 
The complimentary case is studied similarly.\\
Below we assume that $h$ is sufficiently small.
\subsection{Two geometric lemmas}
\label{sec:two-geometric-lemmas}
To prove Theorem~\ref{main:th}, we shall use the following two lemmas.
\begin{Le}
  \label{le:K1}
  There exist two curves in $S_{1,2}$ precanonical with respect to
  $p_2$ and having common endpoints, and $\U_1\subset U$, a 
  neighborhood of $z_0$, such that
  \begin{itemize}
  \item the domain $K_1$ bounded by the two curves is simply
    connected,
  \item $K_1\cap \U_1=S_{1,2}\cap \U_1$.
  \end{itemize}
\end{Le}
\noindent and
\begin{Le}
  \label{le:K0}
  There are exist two curves in $\sigma_2\cup S_{0}\cup \sigma_1$ 
  precanonical with respect to $p_0$ and having common endpoints,
  and $\U_0\subset U$, a  neighborhood of $z_0$, such that
  \begin{itemize}
  \item the domain $K_0$ bounded by the two curves is simply
    connected,
  \item $K_0\cap \U_0= (\sigma_2\cup S_{0}\cup \sigma_1)\cap \U_0$.
\end{itemize}
\end{Le}
\noindent We prove  these two lemmas  in section~\ref{ss:K}.\\
We define $\U=\U_0\cap \U_1$.
\subsubsection{Solution $\psi_1$}
\label{sec:solut-psi_1-psi_2}
We denote by $\psi_{0,0}$ and $\psi_{1,0}$ the solutions $\psi_0$ and
$\psi_1$ constructed by Theorem~\ref{th:CD} for the domain $K_0$, and
consider the solution $\psi_1$ constructed in Theorem~\ref{th:CD} for
the domain $K_1$. In view of Corollary~\ref{cor:th:CD}, in $\U$
(possibly reduced somewhat), all the three solutions are analytic,  
the Wronskian of $\psi_{0,0}$ and $\psi_{1,0}$ does not vanish 
(see also Lemma~\ref{le:2}), and one has
\begin{equation}
  \label{eq:three-psi}
  \psi_1=a\psi_{1,0}+b\psi_{0,0},
\end{equation}
where $a$ and $b$ are $h$-periodic coefficients (see
section~\ref{s:space-of-sol}).  We prove
\begin{Le}
  \label{le:a-and-b} 
  One can reduce $\U$ so that for $z\in \U$
  \begin{equation}
    \label{est:a-and-b}
    a(z)=1+O(h^{L+\frac23}),\quad b(z)=O(h^{L+\frac23}),\qquad h\to 0.
  \end{equation}
\end{Le}
\begin{proof}
  In $\U$ (possibly reduced somewhat), the coefficients $a$ and $b$ are described
  by~(\ref{eq:periodic-coef})  with $\psi=\psi_1$, $f=\psi_{1,0}$ and
  $g=\psi_{0,0}$.\\
  Let $\gamma_{12}= (\sigma_1\cup \sigma_2)\cap \U$.  By
  Lemmas~\ref{le:K1} and~\ref{le:K0} one has $\gamma_{12}\subset K_0$ and
  $\gamma_{12}\subset K_1$. \\
  First, we fix $c\in(1,2)$ and assume that
  $\{z,z+h\}\subset(\gamma_{12})^{ch}$. \\
  In view of Lemma~\ref{le:actions-signs}, one has
  $|\rho_1|=|\rho_2|=1$ on $\gamma_{12}$. This and the definitions of
  $|\rho_1|$ and $|\rho_2|$, see section~\ref{as-sol:est}, imply that
  there exists $C>0$
  such that $|\rho_1(z)|,|\rho_2(z)|\le C$ in $(\gamma_{12})^{ch}$. \\
  As $(\gamma_{12})^{ch}$ is a subset of both $K_0^{ch}$ and
  $K_1^{ch}$, by means of~(\ref{eq:sol-as-exp}) and~(\ref{eq:W-est}),
  we get
  \begin{equation*}
    a=\frac{(\psi_1,\,\psi_{0,0})}{(\psi_{1,0},\psi_{0,0})}=
    \frac{( W_1+O( h^{L+\frac43}),\, W_0+O( h^{L+\frac43}))}
    {( W_1+O( h^{L+\frac43}),\, W_0+O( h^{L+\frac43}))}=
    \frac{( W_1, W_0)+O(h^{L+\frac53})}
    {( W_1, W_0)+O(h^{L+\frac53})}.
  \end{equation*}
  Lemma~\ref{le:2}, then, yields the asymptotic representation for $a$
  from~\eqref{est:a-and-b}. Reasoning similarly, we get
  \begin{equation*}
    b=\frac{(\psi_{1,0},\psi_1)}{(\psi_{1,0},\psi_{0,0})}=
    \frac{( W_1, W_1)+O(h^{L+\frac53})}
    {( W_1, W_0)+O(h^{L+\frac53})}=O(h^{L+\frac23}).
  \end{equation*} 
  This is the estimate for $b$ from~\eqref{est:a-and-b}.
  
  Let $c_1$ and $c_2$ correspond to the minimal strip $S(c_1,c_2)$
  containing $(\gamma_{12})^{ch}$.  We proved estimates~\eqref{est:a-and-b}
  for $a(z)$ and $b(z)$ in the  case where $z,z+h\in (\gamma_{12})^{ch}$.  
  As $c>1$ and as $a$ and $b$ are $h$-periodic, 
  these estimates remain valid in $S(c_1,c_2)$. This implies Lemma~\ref{le:a-and-b}.
\end{proof}
In view of Lemma~\ref{le:K1}, the solution $\psi_1$ admits
representation~\eqref{eq:sol-as-exp} with $l=1$ in $S_{1,2}\cap \U$.
Let us prove that it admits this representation in $S_0\cap\U$.\\
In view of Lemma~\ref{le:K0}, the solutions $\psi_{0,0}$ and
$\psi_{1,0}$ admit representations~\eqref{eq:sol-as-exp} with $l=0$
and $l=1$ in $S_{0}\cap \U$.  Substituting~\eqref{est:a-and-b} and
these representations into~\eqref{eq:three-psi} and
using~(\ref{eq:W-est}), we get for $z\in S_0\cap \U$
\begin{align*}
  \psi_1(z)&=(1+O(h^{L+\frac23}))(W_1(z)+O(h^{L+1+\frac13}\rho_1(z)))\\
           &\hspace{3.5cm}+
             O(h^{L+\frac23})(W_0(z)+O(h^{L+1+\frac13}\rho_0(z)))\\
           &=W_1(z)+O(h^{L+1}\rho_1(z))+O(h^{L+1}\rho_0(z)).
\end{align*}
In view of Lemma~\ref{le:actions-signs}, in $S_0$ one has
$|\rho_0(z)|\le |\rho_1(z)|$. For $\psi_1$ in $S_{0}\cap \U$, this implies
representation~\eqref{eq:sol-as-exp}  with $L$
replaced with $L-1$. As we can increase $L$, we actually
proved~\eqref{eq:sol-as-exp} for $\psi_1$ in the whole domain $\U$.\\
Now, we note that
\begin{equation}\label{help}
  h^{\frac12}|\rho_0(z)|\le C | h^{\frac13}w_0(h^{-\frac23}\zeta(z))|+ 
  C | h^{\frac13}w_0'(h^{-\frac23}\zeta(z))|,\quad z\in U.
\end{equation}
For sufficiently large values of $h^{-\frac23}|\zeta(z)|$, this
estimate follows from the definition of $\rho_0$ and the asymptotic
formulas~(\ref{ai:as:simple}). For bounded $h^{-\frac23}|\zeta(z)|$,
it follows from the fact that $w$ and $w'$ do not have common zeros.\\
Estimates~\eqref{eq:sol-as-exp} and~\eqref{help}
imply~\eqref{eq:asymptotics-main-w:0} with $L$ replaced with $L-1$. As
we can increase $L$, this completes the proof of the statement of
Theorem~\ref{main:th} on the solution $\psi_1$ in the case that we
consider.
\subsubsection{Solution $\psi_0$} 
\label{sec:solution-psi_0}
Let $\psi_{1,1}$ and $\psi_{2,1}$ be the solutions $\psi_1$ and
$\psi_2$ constructed by Theorem~\ref{th:CD} for the domain $K_1$, and let
 $\psi_0$ be the solution constructed by Theorem~\ref{th:CD} for the
domain $K_0$. For $z\in \U$ (possibly reduced somewhat) one has
\begin{equation}
  \label{eq:three-psi-new}
  \psi_0=a\psi_{1,1}+b\psi_{2,1},
\end{equation}
where $a$ and $b$ are $h$-periodic. One proves
\begin{Le}
  \label{le:a-and-b-new} One can reduce $\U$ so that, for $z\in \U$,  one has
  \begin{equation}
    \label{est:a-and-b-new}
    a(z)=-1+O(h^{L+\frac23}),\quad b(z)=-1+O(h^{L+\frac23}),\quad h\to 0.
  \end{equation}
\end{Le}
\begin{proof} We omit details explained in the course of the proof of
  Lemma~\ref{le:a-and-b}. We fix $c\in (1,2)$ and assume that 
  $z,\;z+h\in(\gamma_{12})^{ch}$. 
  For the coefficient $a$ from~(\ref{eq:three-psi-new}), we get
  \begin{equation*}
    a=\frac{(\psi_0,\,\psi_{2,1})}{(\psi_{1,1},\psi_{2,1})}=
    \frac{( W_0,\, W_2)+O(h^{L+\frac53})}
    {( W_1,\, W_2)+O(h^{L+\frac53})}=
    \frac{(-W_1-W_2,\, W_2)+O(h^{L+\frac53})}
    {( W_1,\, W_2)+O(h^{L+\frac53})},
  \end{equation*} 
  where, in the last step, we used relation~(\ref{eq:three-W}).
  Continuing, we get $a=-1+O(h^{L+\frac23})$. Similarly one proves
  that $b=-1+O(h^{L+\frac23})$.  So,~(\ref{est:a-and-b-new}) is proved
  for $z$ we considered. Reasoning as in the completion of the proof
  of Lemma~\ref{le:a-and-b}, we complete the proof of
  Lemma~\ref{le:a-and-b-new}.
\end{proof}
By Theorem~\ref{th:CD} and Lemma~\ref{le:K0}, the solution $\psi_0$ admits
representation~\eqref{eq:sol-as-exp} with $l=0$ in
$(\sigma_1\cup S_{0}\cup \sigma_2)\cap \U$. Estimates~\eqref{est:a-and-b-new}
and~(\ref{eq:W-est}) imply that in $S_{1,2}\cap \U$ one has
\begin{equation*}
  \psi_0=-W_1-W_2+O((|\rho_1|+|\rho_2|)h^{L+1})
  =W_0+O((|\rho_1|+|\rho_2|)h^{L+1}). 
\end{equation*}
In view of Lemma~\ref{le:actions-signs} and the definitions of
$|\rho_j|$, in $S_{1,2}$ one has $|\rho_1|+|\rho_2|\le C|\rho_0|$ which
yields~(\ref{eq:sol-as-exp}) with $l=0$ in $S_{2}\cap \U$. Reasoning
as in the completion of section~\ref{sec:solut-psi_1-psi_2}, we
complete the proof of Theorem~\ref{main:th} for $\psi_0$ in the case
that we consider.
\subsubsection{Solution $\psi_2$}
One proves the main theorem for  $\psi_2$ using the same techniques as for $\psi_0$ and $\psi_1$.
So, we omit the proof and note only that in $S_0$ one represents $\psi_2$ as a linear combination of 
$\psi_{1,0}$ and $\psi_{0,0}$, and computes  the coefficients in this linear combination as  in the 
case of $\psi_0$.
%
\begin{figure}
  \centering
  \includegraphics[height=8cm]{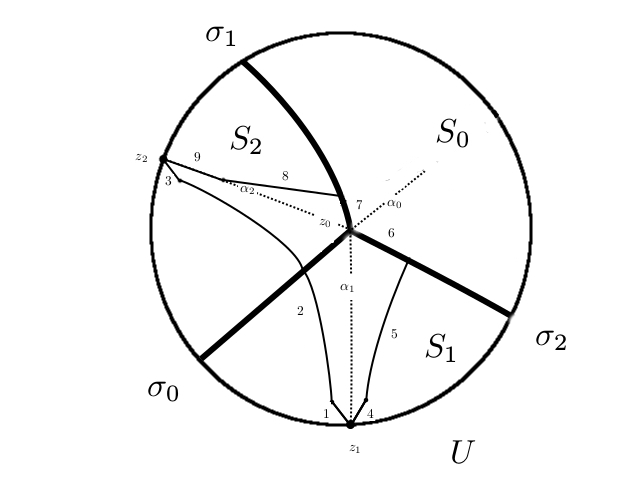}
  \caption{Domain $K_1$}\label{fig:K1}
\end{figure}
%
\section{Proof of the geometric lemmas}
\label{sec:proof-geom-lemm}
\subsection{Proof of Lemma~\ref{le:K1}}
\label{ss:K}
This is done in several steps. \\
Below, all the precanonical lines are precanonical with 
respect to the branch $p_2$ of the complex momentum.
We recall that  $p_2$ is defined and analytic in the domain $U_2$ 
and continuous up to its boundary.
\subsubsection{AntiStokes lines}
\label{sec:antistokes-lines}
We recall that the Stokes lines $\sigma_j$ are defined
by~(\ref{sigma:construction}). The {\it AntiStokes
lines}, $(\alpha_j)_{j\in\Z_3}$, are defined as
\begin{equation}
  \label{alpha:construction}
  \alpha_j:=\zeta^{-1}(V\cap e^{-2\pi i j/3}\,[0,+\infty)).
\end{equation} 
For $j\in\Z_3$, $\sigma_j\cap \alpha_j=\{z_0\}$ and the curve
$\sigma_j\cup \alpha_j$ is analytic. The angles between any two of the
AntiStokes lines at $z_0$ equal $2\pi/3$.\\
In the case we study, the Stokes and AntiStokes lines are
pictured in Fig.~\ref{fig:K1}; the AntiStokes lines are represented by
dotted lines. In particular, $\alpha_2$ goes up from $z_0$,
and $\alpha_1$ goes down from $z_0$.\\
Reducing $U$ if necessary, we assume that the AntiStokes lines 
$\alpha_1$ and $\alpha_2$ are
vertical in $U$. As in Fig.~\ref{fig:K1}, let $z_1$ be
the lower end of $\alpha_1$ and $z_2$ the upper end of $\alpha_2$.\\
 One has
\begin{Le} \label{le:AntiStokes:1} Along the AntiStokes lines
$\alpha_0$, $\alpha_1$ and $\alpha_2$,  one has $\re \int_{z_0}^zp_2\,dz=0$. 
The vector field  $z\mapsto v(z)=\nabla\im \int_{z_0}^zp_2\,dz$ vanishes only at $z=z_0$.  
The AntiStokes lines are tangent to this vector field  at $z\ne z_0$. As
$z$ moves away from $z_0$, $\im \int_{z_0}^zp_2\,dz$ monotonously
increases along $\alpha_2$ and monotonously decreases along $\alpha_1$ and
$\alpha_0$.
\end{Le}
\begin{proof}
  The statement on $\re \int_{z_0}^zp_2\,dz$ 
  follows directly from the definitions of the function $\zeta$ and of 
  the AntiStokes lines. We note that $\|v(z)\|= |p_2(z)|$, and that $p_2(z)$ 
  vanishes only at $z_0$ (the complex momentum 
  vanishes modulo $\pi$ only at turning points and $z_0$ is the only turning point 
  in $U$).  Therefore, the vector field $v$ vanishes only at $z=z_0$.
  The statement on $\re \int_{z_0}^zp_2\,dz$ and the Cauchy-Riemann equations
  imply that  the AntiStokes lines are tangent to  the vector field $v$ where it does 
  not vanish.  This and the first two points of Lemma~\ref{le:actions-signs}
  imply the statements of Lemma~\ref{le:AntiStokes:1} on
  $\im \int_{z_0}^zp_2\,dz$.
\end{proof}
\noindent We also use
\begin{Le}\label{le:alpha}
  There exists $\tilde U\subset U$, a neighborhood of $z_0$, such that the 
  lines $\alpha_1\cap\tilde U$ and $\alpha_2\cap \tilde U$ are precanonical.\\
  Let us parametrize $(\alpha_1\cup\alpha_2)\cap \tilde U$ by \
  $y=\im z$, \ $z=z(y)=x(y)+iy$. Then, if $y\ne \im z_0$, one has
  \begin{gather}
    \label{eq:strictly-can:1}
    \frac{d}{dy}\im \int_{z_0}^{z(y)}p_2(z)\,dz>0,\\
    \label{eq:strictly-can:2} \frac{d}{dy}\im
    \int_{z_0}^{z(y)}(p_2(z)-\pi)\,dz<0.
  \end{gather}
\end{Le}
\begin{proof}  As $\alpha_1$ and $\alpha_2$  are vertical,  
inequality~\eqref{eq:strictly-can:1} follows from Lemma~\ref{le:AntiStokes:1}.
Furthermore, one has
  \begin{equation*}
    \frac{d}{dy}\im \int_{z_0}^{z(y)}(p_2-\pi)\,dz=\im (z'(y)p_2(z))-\pi.
  \end{equation*}
  Therefore, as $p_2(z_0)=0$, reducing $U$ somewhat if necessary, we
  ensure~\eqref{eq:strictly-can:2}.\\
  Since $\alpha_1\cup\alpha_2$ is vertical,~\eqref{eq:strictly-can:1}
  and~\eqref{eq:strictly-can:2} imply that the curve $\alpha_1\cup\alpha_2$ is
  precanonical.
\end{proof}
\noindent Below, we assume that $\tilde U=U$ (if necessary we reduce $U$
somewhat).
\subsubsection{Precanonical line $\gamma_1$}
\label{sec:prec-line-gamm}

We now construct a precanonical line $\gamma_1\subset S_{1,2}$. 
It consists of three segments 1,2 and 3 shown in
Fig.~\ref{fig:K1}. Let us describe them. 
\\
{\it The segments 1 and 3.} \ To construct these segments, we use
\begin{Le}
  \label{le:C1}
  Let $\gamma$ be a compact vertical $C^1$-curve parameterized by
  $y=\im z$, \ $z=z(y)=x(y)+iy$. We assume
  that~\eqref{eq:strictly-can:1}--~\eqref{eq:strictly-can:2} hold
  along $\gamma$. Then, any compact $C^1$-curve sufficiently close in
  $C^{1}$-topology to $\gamma$ is precanonical.
\end{Le}
\noindent This statement immediately follows from the definition of
the precanonical curves.
\\
{\it The segment 1.} It is a segment of a compact precanonical
$C^1$-curve $c_1\subset S_{1,2}$ that begins at $z_1$ and above $z_1$ goes to
the left of $\alpha_1$. When choosing $c_1$, we take  an
internal point of $\alpha_1$ as $\tilde z_1$, and, as $c_1$,  we take a $C^1$-curve close
enough in $C^1$-topology to $\alpha_1$ between $z_1$ and
$\tilde z_1$.  Lemmas~\ref{le:alpha} and~\ref{le:C1} guarantee that
$c_1$ is a precanonical line.\\
{\it The segment 3.} Similarly, the segment 3 is a segment of a
compact precanonical $C^1$-curve $c_3\subset S_{1,2}$, having the upper end
at $z_2$ and going to the left of $\alpha_2$ below the point $z_2$.\\
{\it The segment 2.} We note that $\alpha_1\cup\alpha_2$ is a level curve of the
harmonic function $z\to\re \int_{z_0}^zp_2(z)\,dz$ in $S_{1,2}$. The
segment 2 is a segment of another level curve $c_2$ of this function
in $S_{1,2}$. This curve is located to the left of $\alpha_1\cup\alpha_2$. It
does not contain the point $z_0$, the only point in $S_{1,2}$ where
$p_2$ vanishes. So,  $c_2$ is smooth.
We choose $c_2$ sufficiently close to $\alpha_1\cup\alpha_2$ to ensure that 
\begin{itemize}
\item $c_2$ is vertical (as $\alpha_1$ and $\alpha_2$ are);
\item one has~\eqref{eq:strictly-can:1} along $c_2$ (the vector field
  $\nabla\im \int_{z_0}^zp_2(z)\,dz$ does not vanish along $c_2$ and 
  is tangent to $c_2$);
\item \eqref{eq:strictly-can:2} holds along $c_2$ (as it holds along
  $\alpha_1\cup \alpha_2$);
\item $c_2$ intersects both $c_1$ and $c_3$.
\end{itemize}
Clearly, $c_2$ is precanonical.\\
{\it The curve $\gamma_1$.}   The segment 1 is the 
segment of $c_1$ between $z_1$ and the point
of intersection of $c_1$ and $c_2$, the segment 2 is the segment of
$c_2$ between the segment 1 and the point of intersection of $c_2$ and
$c_3$, and the segment 3 is the segment of $c_3$ connecting the
segment 2 with $z_2$. Clearly, the curve
$\gamma_1$ made of segments 1--3 is precanonical.
\subsubsection{The sign of $\im p_2$ in $S_{2}$}
\label{sec:sign-im-p_0}
The only place where we use our assumption on the direction 
of the tangent vector to $\sigma_0$ at $z_0$ is the proof of
\begin{Le} 
  \label{le:im_p_0}
  Both in $S_{2}$ between the curves $\alpha_2$ and $\sigma_1$ and 
  on these curves, near $z_0$ one has $\im p_2(z)<0$ if $z\ne z_0$.
\end{Le}
\begin{proof}
  Below we assume that either $z$ is in $S_{2}$ between the curves $\alpha_2$ 
  and $\sigma_1$ or on  one of these curves.  In view
  of~(\ref{eq:p-z}), we can write
  \begin{equation}
    \label{eq:p0:loc}
    p_2(z)=k_1\tau(1+O(\tau)), \quad 
    \int_{z_0}^zp_2(z)\,dz=\frac23k_1\tau^3(1+O(\tau)), \quad 
    \quad z\to z_0,
  \end{equation}
  where $k_1\ne 0$ and $\tau$ is the branch of $\sqrt{z-z_0}$ analytic in $U_2$ 
  and positive if $z> z_0$.\\
  Let $0<\theta_2<\pi$ be the angle at $z_0$ between  the line $\{z\ge z_0\}$ 
  and the curve $\alpha_2$ . Note that the angle between $\sigma_0$ and $\alpha_2$ 
  equals $\pi/3$. Therefore, as the tangent vector to $\sigma_0$ at $z_0$ is either 
  directed downwards or parallel to the real line and  directed to the left,  one 
  has  $2\pi/3\le\theta_2<\pi$.\\
  In view of Lemma~\ref{le:AntiStokes:1}, along $\alpha_2$, \ $\re \int_{z_0}^zp_2dz=0$ 
  and  $\im \int_{z_0}^zp_2dz$ is monotonously increasing. This and the second formula 
  in~\eqref{eq:p0:loc} imply that
  \begin{equation}
    \label{arg-k-1}
    \arg k_1+\frac32 \theta_2=\frac\pi2 \mod2\pi.
  \end{equation}
  Let $z-z_0=|z-z_0|e^{i\theta}$. 
  Using~\eqref{arg-k-1} and the first formula in~\eqref{eq:p0:loc}, we
  get near $z_0$
  \begin{equation}
    \label{eq:2}
    \frac{\im p_2(z)}{|p_2(z)|}=
    \sin\left(\arg k_1+\frac\theta2+o(1)\right)
    =\cos\left(\theta_2-\frac{\theta-\theta_2}2+o(1)\right).
  \end{equation}
  Now, we note that, for $z$ we consider, near $z_0$ 
  one has  $\theta_2-\pi/3+o(1)\le \theta\le \theta_2+o(1)$. 
  Therefore, for $z$ sufficiently close to  $z_0$, one has 
  \begin{equation*}
    \frac{2\pi}3+o(1)\le \theta_2+o(1)
    \le\theta_2-\frac{\theta-\theta_2}2 \le
    \theta_2+\frac\pi6+o(1)<\frac{7\pi}6.   
  \end{equation*}
  This and~\eqref{eq:2} implies the statement of
  Lemma~\ref{le:im_p_0}.
\end{proof}
\subsubsection{Precanonical line $\gamma_2$}
\label{sec:prec-line-gamm-1}
The precanonical line $\gamma_2$ is located in $S_{1,2}$ and consists
of six segments 4--9 shown in Fig.~\ref{fig:K1}. Let us describe
them.
\\
{\it The segments 4-5-6-7.} The segment 4 is a segment of a compact
precanonical $C^1$-curve $c_4\subset S_{1,2}$. This curve begins at $z_1$ and
above $z_1$ goes to the right of $\alpha_1$. It is constructed as the curve
$c_1$ containing the segment 1.\\
The segment 5 is a segment of a level curve $c_5$ of the function
$z\to\re \int_{z_0}^zp_2(z)\,dz$ in $S_{1,2}$. The construction of
$c_5$ is similar to one of $c_2$. The curve $c_5$ is located to the
right of $\alpha_1$. We choose $c_5$ sufficiently close to $\alpha_1$. Then,
$c_5$ is a precanonical curve and intersects both $c_4$ and the Stokes
line $\sigma_2$.\\
The segment 4 is the segment of $c_4$ between $z_1$ and the point of
intersection of $c_4$ and $c_5$. The segment 5 connects this point with
a point of $\sigma_2$.\\
We prove
\begin{Le}
  \label{le:Stokes}
  Let $\gamma$ be a vertical curve, let $a\in\gamma$ and let $p$ be a
  branch of the complex momentum continuous on $\gamma$. If, on
  $\gamma$, either $\im\int_a^zp(z)\,dz=0$ or
  $\im\int_a^z(p(z)-\pi)\,dz=0$, then $\gamma$ is precanonical  
  with respect to $p$.
\end{Le}
\begin{proof}
  Assume that $\im\int_a^zp(z)\, dz=0$ on $\gamma$. Then,
  $z\mapsto\im\int_a^z(p(z)-\pi)\, dz=-\pi\im (z-a)$ is decreasing
  along $\gamma$ when $\im z$ increases. Thus, $\gamma$ is precanonical.\\
  If $\im\int_a^z(p(z)-\pi)\,dz=0$, then
  $z\mapsto\im\int_a^zp(z)\,dz=\im\int_a^z\pi\,dz =\im (z-a)$ is
  increasing along $\gamma$ when $\im z$ increases. Thus, $\gamma$ is precanonical.
\end{proof}
\noindent The segment 6 is the segment of $c_6=\sigma_2$ between the
upper end of the segment 5 and the point $z_0$. The segment 7 is the
segment of $c_7=\sigma_1$ between $z_0$ and an internal point $a$ of
$\sigma_1$.  We describe this point later. Lemma~\ref{le:Stokes}
implies that the segments 6 and 7 are precanonical.
\\
{\it Segment 8.} This segment is a segment of $c_8$, the level curve
$\gamma(a)$ of the harmonic function
$z\to \im\int_{z_0}^z(p_2(z)-\pi)\,dz$ that contains $a \in\sigma_1$. To
choose the segment 8, we check
\begin{Le}
  \label{le:6}
  If $a\in\sigma_1\setminus\{z_0\}$ is sufficiently close to $z_0$,
  then $\gamma(a)$  intersects $\sigma_1$ transversally at $a$, 
  enters at $a$ in $S_2$ going upwards, intersects $\alpha_2$ and, 
  up to  intersection and at the intersection point, remains vertical. 
\end{Le}
\begin{proof}[Proof of Lemma~\ref{le:6}]
  Below, we identify the vectors on $\R^2$ with the complex numbers 
  in the standard way, and the bar denotes complex conjugation.  The 
  Stokes line $\sigma_1$  is tangent to the vector field
  $z\mapsto v_0(z)=\overline{p_2(z)}$ at $z\ne z_0$ ($p_2(z_0)=0$).  The
  curve $\gamma(a)$ is tangent to the vector field
  $z\mapsto v_\pi(z)=\overline{p_2(z)}-\pi$.\\
  Let $a\in\sigma_1\setminus {z_0}$ be sufficiently close to the point $z_0$. 
  In view of Lemma~\ref{le:im_p_0}, $\im p_2(a)< 0$.
  Therefore,  $\gamma(a)$ is vertical at $a$. Moreover,
  both the vectors $v_0(a)$ and $v_\pi(a)$ are directed upwards and
  $v_\pi(a)$ is directed to the left of $v_0(a)$. Therefore, at $a$,
  the curve $\gamma(a)$ intersects $\sigma_1$ transversally and enters
  $S_2$ going upwards.\\
  Furthemore, in view of Lemma~\ref{le:im_p_0}, as long as $\gamma(a)$ 
  stays in $S_2$ near $z_0$ between the curves $\alpha_2$ and $\sigma_1$ 
  or on them, it remains vertical.\\
  To complete the proof, it suffices to show that if $a$ is
  sufficiently close to $z_0$, then $\gamma(a)$ intersects $\alpha_2$
  remaining vertical. Therefore, we note that 
  $v_\pi(z_0)=-\pi$.  So, at $z_0$ the vector tangent to $\gamma(z_0)$ 
  is parallel to $\R$, and the curve $\gamma(z_0)$  intersects the analytic 
  curve $\alpha_2\cup\sigma_2$ transversally. Depending  continuously 
  on $a$,  \ $\gamma(a)$ intersects this curve also  for all $a$ sufficiently 
  close to $z_0$. But,  if  $\im a>\im z_0$ and $a$ is sufficiently close to  
  $z_0$, the curve $\gamma(a)$ goes upward from $a$.  Therefore, for $a$
  sufficiently close to $z_0$, the curve $\gamma(a)$ intersects $\alpha_2$
  still remaining vertical. This completes the proof of 
  Lemma~\ref{le:6}.
\end{proof}
\noindent {\it The segments 8 and 9.} 
We choose the point $a$, the end of the segment 7 and the beginning of
the segment 8, so that $c_8=\gamma(a)$ intersects $\alpha_2$ as described
in Lemma~\ref{le:6}. The end of the segment 8 is the point of
intersection of $c_8$ and $\alpha_2$. By Lemma~\ref{le:Stokes}, the segment
8 is precanonical. The segment 9 is the segment of $\alpha_2$ connecting
the upper end of the segment 8 to the point $z_2$. It precanonical
by Lemma~\ref{le:alpha}.

{\it The domain $K_1$} bounded by  $\gamma_1$ and $\gamma_2$
is the one described in Lemma~\ref{le:K1}, the proof of which is
complete.
\subsection{Proof of Lemma~\ref{le:K0}}\label{ss:K3}
\begin{figure}
\centering
\includegraphics[height=7.3cm]{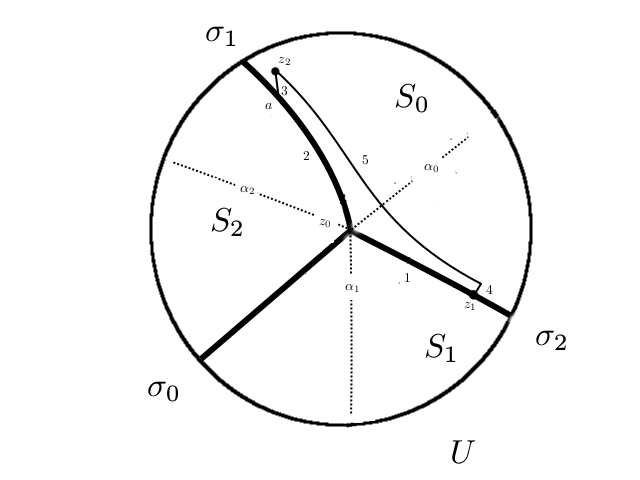}
\caption{Domain $K_0$}\label{fig:K0}
\end{figure}

The proof uses the same techniques as the proof of Lemma~\ref{le:K1}.
Therefore, we omit standard details.
The construction of the curves  $\gamma_1$ and $\gamma_2$
bounding the domain $K_0$ from Lemma~\ref{le:K0} is illustrated 
by Fig.~\ref{fig:K0}. Below, all the precanonical lines are precanonical 
with respect to $p_0$.
%
%
\subsubsection{Curve $\gamma_1$}\label{sss:K0:gamma1}
This  curve consists of  segments 1--3. Let us describe them.\\
We take  an internal point of $\sigma_2$  as $z_1$, and  we fix
 $a$, an internal point of $\sigma_1$. 
The segment 1 is the segment of $\sigma_2$
between $z_1$ and $z_0$, and the segment 2 is the segment of $\sigma_1$ 
between $z_0$ and $a$. \\
To describe the segment 3, we consider  $\gamma_0(a)$, 
the curve in $\sigma_2\cup S_0\cup \sigma_1$ described by the equation  
$\im \int_{a}^z(p_0(z)-\pi)\,dz=0$. We suppose that  $a$ is sufficiently close 
to $z_0$. Then, $\gamma_0(a)$ intersects $\sigma_1$ at $a$ transversally,
enters in $S_0$ going upwards and is vertical 
in  a neighborhood $a$ (To prove this, one uses the observation that near 
$z_0$ on $\sigma_1$ one has $\im p_0(z)>0$. The proof of this observation  
is similar to one of Lemma~\ref{le:im_p_0}.)
The segment 3 is a segment of $\gamma_0(a)$ connecting 
in this neighborhood $a$ to a point $z_2\in S_2$ . We choose $z_2$ later. \\
Lemma~\ref{le:Stokes} imply that the segments 1--3 are precanonical.\\
The points $z_1$ and $z_2$ are the ends of $\gamma_1$.
\subsubsection{Curve $\gamma_2$}
This curve consists of two segments, segments 4  and 5.
\\
The segment 4 is a  segment  of $c_4$, a level curve of the function  
$z\to \re\int_{z_1}^zp_0(z)\,dz$ in $S_0\cup \sigma_2$ that contains the 
point $z_1$. The curve $c_4$ is orthogonal  to $\sigma_2$  at $z_1$. 
\\
Let us note that, under our assumptions on $\sigma_0$ and $\sigma_2$
(see the very beginning of section~\ref{s:proof:main}),  the angle at $z_0$
between $\sigma_2$ and the horizontal line $\{z\ge z_0\}$ belongs to $(0, \pi/3)$.
Possibly reducing $U$ somewhat, we assume that, at any point $\zeta\in\sigma_2$,
the angle between $\sigma_2$ and the line $\{z\ge \zeta\}$ belongs to $(0,\pi/2)$.
Then, $c_4$ is vertical   at least in a neighborhood of the point $z_1$ and 
goes upward from $z_1$ into $S_0$. 
\\
The segment 5 is a segment of a level curve $c_5$ of the function 
$z\to \im\int_{z_0}^zp_0(z)\,dz$ in $S_0$. It is located to the right of 
$\sigma_2\cup\sigma_1$  (which is also a level curve of this function).
We  choose the curve $c_5$ sufficiently close to $\sigma_2\cup\sigma_1$.
Then it is vertical, intersects $\gamma_0(a)$ and $c_4$, and the segments
of these curves between $\sigma_2\cup\sigma_1$ and the intersection points 
are vertical. 
\\
The point $z_2$ is the point of intersection
of $\gamma_0(a)$ and $c_5$. The segment 4 is the segment of $c_4$ between 
$\sigma_2\cup\sigma_1$ and $c_5$, and the segment 5 is the segment of $c_5$
connecting $c_4$ to $z_2$.
\\
The segment 5 is precanonical in view of Lemma~\ref{le:Stokes}.
Arguing as when proving Lemma~\ref{le:alpha} and reducing somewhat $U$ if necessary,
we check that the segment 4 is precanonical.

{\it The domain $K_0$} bounded by the curves 
$\gamma_1$ and $\gamma_2$, is the one
described in Lemma~\ref{le:K0}.  Its proof is complete.
\bibliographystyle{plain} 
\end{document}